\documentclass[a4wide]{article}
\usepackage[utf8]{inputenc}
\usepackage[a4paper, total={6.25in, 9.25in}]{geometry}

\usepackage{hyperref}       
\usepackage{url}            
\usepackage{booktabs}       
\usepackage{amsfonts}       
\usepackage{nicefrac}       
\usepackage{microtype}      
\usepackage{xcolor}         
\usepackage{amsmath}
\usepackage{graphicx}
\usepackage{comment}
\usepackage{algorithm}
\usepackage{algpseudocode}
\usepackage[auth-lg]{authblk}

\newtheorem{remark}{Remark}

\usepackage[page]{appendix} 

\begin{document}

\title{Long-time prediction of nonlinear parametrized dynamical systems by deep learning-based reduced order models}

%

\author[1]{Federico Fatone}
\author[1]{Stefania Fresca}
\author[1]{Andrea Manzoni}

\affil[1]{MOX - Dept. of Mathematics, Politecnico di Milano, P.za Leonardo da Vinci 32, 20133 Milano, Italy\footnote{\texttt{{federico.fatone,stefania.fresca,andrea1.manzoni}@polimi.it}}}

 \date{}\renewcommand\Affilfont{\small}

\maketitle

\begin{abstract}
 Deep learning-based reduced order models (DL-ROMs) have been recently proposed to overcome common limitations shared by conventional ROMs -- built, \textit{e.g.}, exclusively through proper orthogonal decomposition (POD) -- when applied to nonlinear time-dependent parametrized PDEs. In particular, POD-DL-ROMs can achieve extreme efficiency in the training stage and faster than real-time performances at testing, thanks to a prior dimensionality reduction through POD and a DL-based prediction framework. Nonetheless, they share with conventional ROMs poor performances regarding time extrapolation tasks. This work aims at taking a further step towards the use of DL algorithms for the efficient numerical approximation of parametrized PDEs by introducing the $\mu t$-POD-LSTM-ROM framework. This novel technique extends the POD-DL-ROM framework by adding a two-fold architecture taking advantage of long short-term memory (LSTM) cells, ultimately allowing long-term prediction of complex systems’ evolution, with respect to the training window, for unseen input parameter values. Numerical results show that this recurrent architecture enables the extrapolation for time windows up to 15 times larger than the training time domain, and achieves better testing time performances with respect to the already lightning-fast POD-DL-ROMs. 
\end{abstract}

\vspace{1cm}

\section{Introduction}
\label{sec:introduction}

Parameterized PDEs are extensively used for the mathematical description of several physical phenomena. Some instances include fluid dynamics, heat transfer, waves and signal propagation and interference, structure dynamics (including microsystems) and chemically reacting flows \cite{reviewPBM,quarteroni2016reduced,benner2017model}. 
However, traditional high-fidelity, full order models (FOMs) employed for their numerical approximation, such as those based on the finite element method, become infeasible when dealing with complex systems and multiple input-output responses need to be evaluated (like, \textit{e.g.}, for uncertainty quantification, control and optimization) or real-time performances must be achieved. In fact, despite being accurate up to a desired tolerance, they entail unaffordable computational times, sometimes even orders of magnitude higher than the ones required by real-time computing \cite{quarteroni1994numerical}. 

In this context, projection-based reduced order models (ROMs) -- such as POD-Galerkin ROMs -- have been introduced with the goal of enhancing efficiency in timing critical applications. ROMs rely on a suitable offline-online computational splitting, aimed at moving in the offline stage, \textit{i.e.}, the phase where the model is trained and refined before its deployment, the computationally expensive tasks in order to make the online one, \textit{i.e.}, the phase where the model is used for the solution of the problem for new parametric instances, extremely efficient. These methods rely on the assumption that the parameterized PDE solutions manifold can be represented by the span of a small number of basis functions built starting from a set of FOM solutions (computed in the offline stage), the so-called \textit{reduced manifold}. This pipeline allows for a significant dimensionality reduction of the PDE problem and a consequent speed-up in its numerical solution time, for example a $O(10^2)$ speed-up is achieved for Navier-Stokes application \cite{manzonispeedup,ballarin15,dalsanto2019}, and even more for applications in structural mechanics \cite{FGMQ_19,frangi2022}.  

Nevertheless, despite being physics-driven, POD-Galerkin ROMs show severe limitations when addressing nonlinear time-dependent PDEs, which might be related to {\em (i)} the need to rely on high-dimensional linear approximating trial manifolds, {\em (ii)} the need to perform expensive hyper-reduction strategies, or {\em (iii)} the intrinsic difficulty to handle compex physical patterns with a linear superimposition of modes \cite{quarteroni2016reduced}. Furthermore, usually such ROMs do not allow for an effective extrapolation in time, requiring extremely long offline stages in order to compute FOM snapshots defined on a sufficiently long time domain. 

To overcome these drawbacks, several nonlinear -- and in particular artificial neural networks (ANNs) based -- methods have been massively considered  to provide fast approximation of PDE solutions in the last few years, and even before. For instance, the possibility to approximate differential equations solutions through ANNs had already been proposed in \cite{lagaris98} and \cite{Aarts2001}, relying on the universal approximation theorem \cite{Hornik91}. Interest on the topic and practical applications increased recently at a fast pace \cite{khoo_lu_ying_2021,michoski20,berg19}, while variations aimed at introducing physics related losses to link more deeply the ANN framework with the underlying physical model, \textit{e.g.}, with the concept of physics-informed neural networks (PINNs) \cite{raissi19,raissi2018deep}, show the significance ANNs are gaining  in scientific computing. Furthermore, some relevant theoretical results concerning the complexity bounds of the problem and the error of the approximation have also been investigated, for example in \cite{opspetsch20,Kutyniok2021,yarotsky17,franco2021}, thus providing a rigorous mathematical framework to the related problems. \\ 

The combination of the extremely accurate approximation capabilities of ANNs \cite{opspetsch20,yarotsky17} and the concept of reduced order modeling led to the introduction of ANN-based ROMs. In particular, the idea consists in using deep learning (DL) algorithms to perform a nonlinear projection onto a suitable reduced order manifold. For instance, in \cite{NN4,NN5,hesthaven18} a DL-based regressor is employed but a linear reduced manifold is still considered. In \cite{osti20} an ANN-inferred correction term is used to increase the accuracy of the linear projection, while \cite{NN1,NN2} approximate the reduced manifold by means of an ANN. In \cite{NN3}, a convolutional autoencoder is considered to model the reduced order manifold, but the advancement in time is performed by means of a quasi-Newton method requiring the approximation of a Jacobian matrix at every time step. 

A recently proposed strategy \cite{fresca2020comprehensive,franco2021} aims at constructing DL-based ROMs (DL-ROMs) for nonlinear time-dependent parametrized PDEs in a non-intrusive way, approximating the PDE solution manifold by means of a low-dimensional, nonlinear trial manifold, and the nonlinear dynamics of the generalized coordinates on such reduced trial manifold, as a function of the time coordinate and the parameters. The former is learnt by means of the decoder function of a convolutional autoencoder (CAE) neural network; the latter through a (deep) feedforward neural network (DFNN), and the encoder function of the CAE. DL-ROMs outperform POD-based ROMs such as the reduced basis method – regarding both numerical accuracy and computational efficiency at testing stage. With the same spirit, POD-DL-ROMs \cite{fresca2021poddlrom} enable a more efficient training stage and the use of much larger FOM dimensions, without affecting network complexity, thanks to a prior dimensionality reduction of FOM snapshots through randomized POD (rPOD) \cite{halko2011finding}, and a multi-fidelity pretraining stage, where different models (exploiting,  e.g., coarser discretizations or simplified physical models) can be combined to iteratively initialize  network parameters. This latter strategy has proven to be effective for instance in the real-time approximation of cardiac electrophysiology problems \cite{fresca20cardiac,frescafrontiers} and problems in fluid dynamics \cite{fmfluids1}. \\

{This work extends the POD-DL-ROM framework \cite{fresca2021poddlrom} in two directions: first, it replaces the CAE architecture of POD-DL-ROM with a long short-term memory (LSTM) based autoencoder \cite{hochreiter1997long,forgetgate}, in order to better take into account time evolution when dealing with nonlinear unsteady parametrized PDEs ($\mu$-POD-LSTM-ROM); second, it aims at performing extrapolation forward in time (compared to the training time window) of the PDE solution, for unseen values of the input parameters -- a task often missed by traditional projection-based ROMs. Our final goal is to predict the PDE solution on a larger time domain $(T_{in}, T_{end})$ than the one, $(0,T)$, used for the ROM training -- here  $0\leq T_{in} \leq T_{end}$ and $T_{end} > T$. To this aim, we train a $\mu$-POD-LSTM-ROM using $N_t$ time instances and approximate the solution up to $N_t + M$ time steps from the starting point, taking advantage of a time series LSTM-based architecture ($t$-POD-LSTM-ROM) besides the $\mu$-POD-LSTM-ROM introduced before. These architectures mimic the behavior of numerical solvers as they build predictions for future times based on the past. Besides this, the implications of the novelties proposed by the present work are multiple. 

In particular, the main advantages concern:
\begin{itemize}

    \item the possible long-term time extrapolation capabilities of the proposed framework, allowing for a faster offline stage, as FOM snapshots defined on a shorter time domain are required to train the model;

    \item the possibility to predict entire sequences instead of single outputs, that makes the presented method even more efficient than the (already faster than real-time) POD-DL-ROM \cite{fresca2021poddlrom},
    
\end{itemize}

\noindent while at the same time preserving the main strengths of POD-DL-ROM \cite{fresca2021poddlrom}, which are:

\begin{itemize}
    \item the possibility to query the method at a specific time for what concerns unsteady dynamical systems, without requiring the computation of the solution at previous time steps, as a traditional time marching method would do;
    
    \item the possibility of using coarser temporal discretizations with respect to the ones used to ensure stability for high-fidelity numerical solvers \cite{fresca20cardiac};

    \item the avoidance of using expensive hyper-reduction techniques often required by POD-based ROMs;
    
    \item the possibility to return outputs depending on selected  problem state variables, without requiring to approximate all of them. \end{itemize}

}

The paper is divided in five sections. In Section \ref{sec:time_ext} we describe the $\mu t$-POD-LSTM-ROM framework used to predict PDE solutions for unseen parameter instances and times. Section \ref{sec:mu_PODLSTMROM} and Section \ref{sec:t_PODLSTMROM} introduce the $\mu$-POD-LSTM-ROM and $t$-POD-LSTM-ROM architectures, the former enriching POD-DL-ROM with LSTM-based autoencoder and the latter providing time extrapolation capabilities to the framework. In Section \ref{sec:results} we report the accuracy results and performances assessments of $\mu t$-POD-LSTM-ROM on three parametrized test cases, namely: {\em (i)} 3 species Lotka-Volterra equations, {\em (ii)} unsteady advection-diffusion-reaction equation, {\em (iii)} incompressible Navier-Stokes equations.


\section{Achieving time extrapolation capabilities with LSTM cells}
\label{sec:time_ext}

After recalling the formulation of a POD-DL-ROM, in this section we address the construction of the proposed $\mu t$-POD-LSTM-ROM framework to predict PDE solutions for unseen parameter instances and times; the main ingredients to reach this goal -- the $\mu$-POD-LSTM-ROM and the $t$-POD-LSTM-ROM architectures -- will instead be detailed in the following sections. 

The space and time discretization on a nonlinear, time-dependent, parametrized PDE problem -- performed, \textit{e.g.}, through a finite element method -- produces a (high-dimensional) dynamical system of the form:
\begin{equation}
\label{eq:FOM}
\begin{cases}
{\bf M}(\boldsymbol{\mu}) \mathbf{\dot{u}}_h(t;\boldsymbol{\mu}) = \mathbf{f}(t, \mathbf{u}_h(t;\boldsymbol{\mu}); \boldsymbol{\mu}), \qquad t \in (0, T),\\
\mathbf{u}_h(0;\boldsymbol{\mu})=\mathbf{u}_0(\boldsymbol{\mu}). 
\end{cases}
\end{equation}
where $\mathbf{u}_h:(0,T) \times \mathcal{P} \rightarrow \mathbb{R}^{N_h}$ is the parametrized solution of (\ref{eq:FOM}), $\mathbf{u}_0 : \mathcal{P} \rightarrow \mathbb{R}^{N_h}$ is the initial datum, $\mathbf{f} : (0,T) \times  \mathbb{R}^{N_h}  \times \mathcal{P} \rightarrow \mathbb{R}^{N_h}$ is a (nonlinear) function, representing the system dynamics and ${\bf M}(\boldsymbol{\mu}) \in \mathbb{R}^{N_h \times N_h}$ is the mass matrix of this parametric FOM, assumed here to be a symmetric positive definite matrix. Here we have denoted by  $N_h$  the dimension of the FOM and by  $\mathcal{P} \subset \mathbb{R}^{n_{\mu}}$ the parameters' space. 

POD-based ROMs exploit singular value decomposition (SVD) of the FOM snapshot matrix $\mathbf{S}$, \textit{i.e.}, the data structure containing the full order solutions (\textit{snapshots}) used for training, in order to build a $N$-dimensional space basis and project the system (\ref{eq:FOM}) on it. In this way, an $N$-dimensional reduced solution manifold $\mathcal{M}_N$ is obtained \cite{quarteroni2016reduced}. Since the solution of the FOM  \eqref{eq:FOM} can be computationally unaffordable, we aim at replacing it with the approximation obtained through suitable ROMs.

Since  \eqref{eq:FOM} entail a temporal evolution of the described phenomenon, the ROMs used to approximate its dynamics must include some kind of time parameter (even dimensionless), or at least some procedure allowing the advancement in time of the solution to work properly. POD-based ROMs, for instance, usually consider a time marching scheme to handle the dynamical system they entail, similarly to the ones used at the FOM level (\textit{e.g.}, finite differences or backward differentiation formulas). POD-DL-ROMs as described in \cite{fresca2021poddlrom} treat instead the time like an additional physical parameter to be provided as input to the feedforward neural network $\boldsymbol{\phi}_n^{FFNN}(\cdot; \cdot, \boldsymbol{\theta}_{FFNN})$ -- being $\boldsymbol{\theta}_{FFNN}$ its trainable parameters -- mapping the parameter vector $\boldsymbol{\mu} \in \mathcal{P} \subset \mathbb{R}^{n_{\mu}}$ to the low ($n$-)dimensional nonlinear manifold, where $n$ is very close or even equal to the intrinsic dimension of the problem ($n_{\boldsymbol{\mu}} + 1$):
\begin{equation*}
\label{eq:with_time}
\mathbf{u}_n(t; \boldsymbol \mu) = \boldsymbol{\phi}_n^{FFNN}(t; \boldsymbol \mu, \boldsymbol{\theta}_{FFNN}),
\end{equation*}
so that in the end, the network learns a mapping of the form
\begin{equation*}
\boldsymbol{\phi}_n^{FFNN}(\cdot ;  \cdot, \boldsymbol{\theta}_{FFNN}) : (0, T) \times \mathbb{R}^{n_{\mu}} \rightarrow \mathbb{R}^n.
\end{equation*}
This latter, starting from each pair $(time, parameters)$, produces the low-dimensional representation of the solution for those particular instances with a direct input-output relation for each time instant.

Despite being fast and accurate, this approach neglects the correlation between consequent time steps of the solution, leaving the opportunity to further increase its efficiency. In fact, the simple selection of an initial condition and of an initial time should ideally contain enough information to reconstruct the entire temporal evolution of the solution. In this context, it would be ideally possible -- and desirable -- to obtain a map under the form 
\begin{equation*}
\boldsymbol{\Lambda}_N(\cdot, \boldsymbol{\theta}_{\Lambda}) : \mathbb{R}^{n_{\mu}} \rightarrow \mathbb{R}^N \times (0, T)
\end{equation*}
that, considering the initial time as $t=0$, would provide the solution for a time horizon $(0,T)$ as long as necessary, thus enhancing time extrapolation capabilities. In this case, time would be considered implicitly by the model; this latter shall then learn the evolution of the problem through its trainable parameters. This would entail the setting of an algorithm -- conceptually closer to a classical numerical solver than what presented before -- as the obtained solution would be a sequence of vectors representing the evolution of the system in time, rather than a single result of a specific time query.

To better fit the working mechanism of numerical solvers, the application of recurrence strategies to ANNs emerges as a suitable solution, as they add to traditional feedforward architectures feedback connections allowing to treat inputs and outputs in the form of sequences. This is expected to enhance the reconstruction of the underlying dynamics, as this latter can be learned by the network implicitly. 

Moreover, such architectures have proven to be effective in time series prediction problems \cite{LSTMtimeseries1} -- even in the context of PDEs approximation \cite{maulik20} -- opening the possibility for ROMs to advance in time with respect to the FOM snapshots they are trained with, performing extrapolation in time. Finally, recurrence mechanisms such as LSTM are also suitable for increasing speed performances at prediction time, as entire long temporal sequences can be returned as output from the architecture, requiring less neural network queries, and thus improving the overall efficiency of the method.

\subsection{Time extrapolation problem}

In the context of ROMs, he problem of time extrapolation hence requires the training of a ROM (eventually based on deep learning) on a training set including snapshots 
\begin{equation}
    \mathbf{u}_h(t; \boldsymbol{\mu}) \ \ \ \ \textnormal{with} \ \  \boldsymbol{\mu} \in \mathcal{P}_{train} \ \textnormal{and} \  t \in (0,T)
\end{equation}
to be used to predict solutions defined on a larger temporal domain $\mathcal{T} = (T_{in}, T_{fin})$ with $0 \le T_{in} \le T_{fin}$ and $T_{fin} > T$:
\begin{equation}
    \mathbf{u}_h(t; \boldsymbol{\mu}) \ \ \ \ \textnormal{with} \ \  \boldsymbol{\mu} \in \mathcal{P}_{test} \ \textnormal{and} \  t \in \mathcal{T}.
\end{equation}
Parameters $\mathcal{P}_{test} \subset \mathcal{P}$ and $\mathcal{P}_{train} \subset \mathcal{P}$ are such that $\mu_i^{train, min} \le \mu_i^{test} \le \mu_i^{train, max} \ \ \forall i \in \{1, \dots, n_{\boldsymbol{\mu}}\}$.\\

Also for this case, an analogy with time marching numerical solvers can be found. Indeed, traditional numerical solvers integrate in time the system of PDEs  starting from the given initial condition $\mathbf{u}_h(0; \boldsymbol{\mu})$, and build the solution iteratively exploiting past solution's values found by the solver itself. Our goal is to build a DL-based solution {\em approximator} able to proceed in time in a similar fashion. Ultimately, our objective is to train the framework using (the first) $N_t$ time steps of the FOM solution, to get the solution of the problem up to $N_t + M$ time steps from the starting point.

\subsection{$\mu t$-POD-LSTM-ROM architecture}
\begin{figure}[b!]
    \centering
        \vspace{-0.45cm}
  \includegraphics[width=0.85\textwidth]{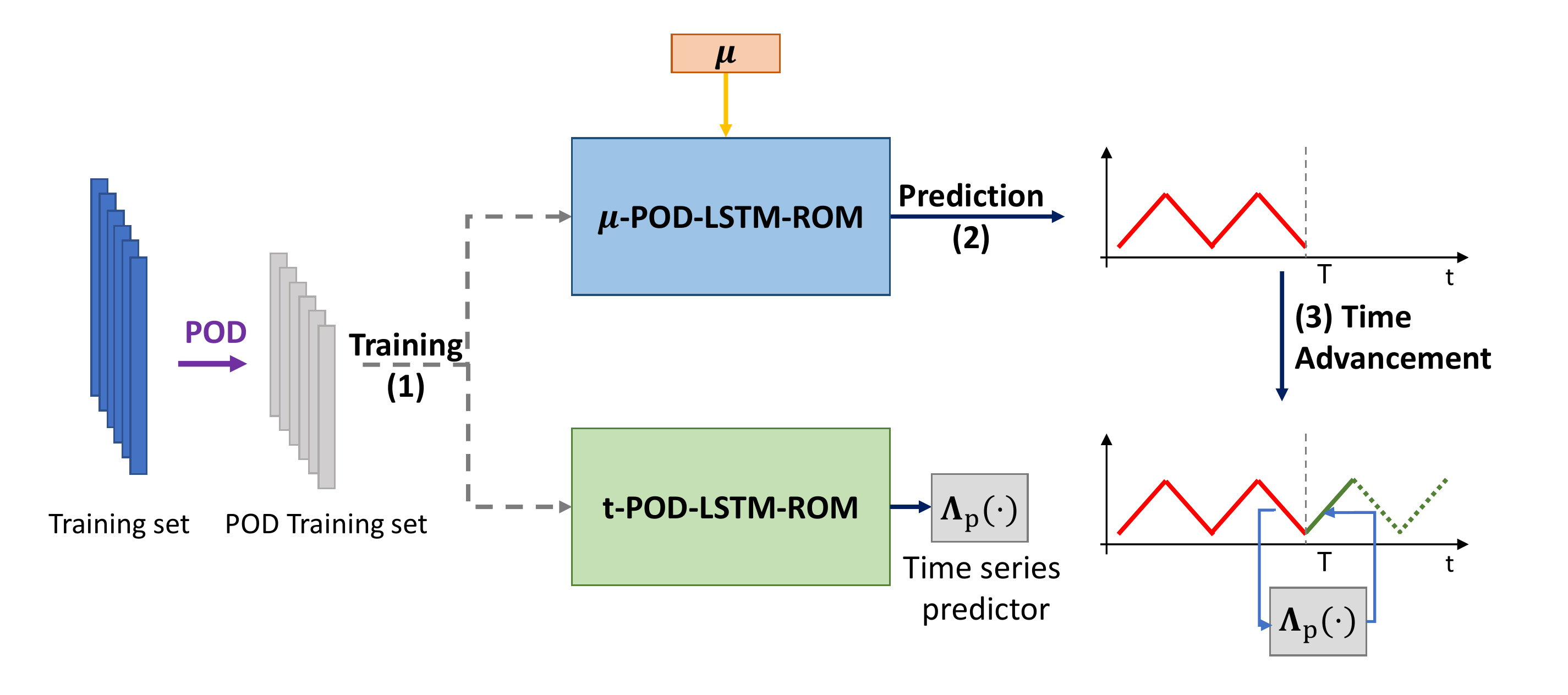}
    \vspace{-0.45cm}
    \caption{The $\mu t$-POD-LSTM-ROM framework. \textit{(1) Training:} both $\mu$-POD-LSTM-ROM and $t$-POD-LSTM-ROM are trained on the same set of FOM snapshots reduced by means of (r)POD; \textit{(2) Prediction:} The $\mu$-POD-LSTM-ROM is employed to predict the (r)POD coordinates on the time interval $(0,T)$ on which the training snapshots were defined for new parameters instances; \textit{(3) Time Advancement:} starting from the sequence predicted by the $\mu$-POD-LSTM-ROM, the $t$-POD-LSTM-ROM time series predictor is used to advance in time and perform time extrapolation.}
    \label{cpodlstmromfig}
\end{figure}

The problem addressed in this paper is therefore two-fold: it deals with {\em (i)} the prediction of the solution of the parametric PDE problem for a new instance of the parameters' space belonging to the set $\mathcal{P}_{test} \subset \mathcal{P} \in \mathbb{R}^{n_{\mu}}$ (solution inference for new parameters values) and {\em (ii)} the forecast of the temporal evolution of that solution for unseen times (time extrapolation).
A natural way to tackle this problem is to pursue a {\em divide-and-conquer} strategy, splitting its solution into a two-steps process that exploits two paired LSTM-based ANN architectures: the first one ($\mu$-POD-LSTM-ROM) addressing the issue of predicting the solution for unseen parameters;  the second one ($t$-POD-LSTM-ROM, where "t" stands for \textit{time series})  extending the solution in time, starting from the sequence predicted by the former. The resulting technique ($\mu t$-POD-LSTM-ROM), summarized in Figure \ref{cpodlstmromfig}, can be described as follows:
\begin{itemize}
    \item The training stage is performed in parallel for the two paired architectures on the same dataset obtained from FOM solutions after a first POD-based dimensionality reduction. In the end, $\mu$-POD-LSTM-ROM will produce a structure able to predict the solution for unseen parameters, but on the same time domain considered during training, while $t$-POD-LSTM-ROM will produce a time series predictor $\boldsymbol{\Lambda}_p(\cdot)$ that takes $p$ time steps from the past and $\boldsymbol{\mu}$, and returns as output the forecast for $k$ time steps in the future for that $\boldsymbol{\mu}$ value. From now on, the two architectures will act as separate entities;
    
    \item $\mu$-POD-LSTM-ROM takes the vector $(t_i, \boldsymbol{\mu})$ of a starting time in the interval $(0, T)$ (discretized in $\{t_0, \dots, t_{N_t-1}\}$) and of the parameters' instance, and performs the approximation of the solution on this time domain seen during the training stage. This step produces very accurate outputs that are also particularly smoother w.r.t. POD-DL-ROM ones thanks to the LSTM architecture producing sequences as outputs. This enhances the performances of the time series predictor, as it would potentially incur in stability issues by propagating the small oscillations somehow unavoidable in the POD-DL-ROM framework without LSTM cells. Note that all the predictions at this stage are performed in the reduced dimension $N$ (the POD basis one);
    
    \item $t$-POD-LSTM-ROM takes the last $p$ time steps of the $\mu$-POD-LSTM-ROM predicted sequence and predicts the $k$ following ones. Then, it takes the last $p$ time steps of the new predicted sequence and advances of other $k$ steps, and keeps advancing the prediction in this way. This strategy thus performs extrapolation in time, virtually with no final time limit, acting as an auto-regressive model. Note that the $t$-POD-LSTM-ROM architecture is general enough to be used also on top of other ROMs -- \textit{e.g.}, POD-DL-ROMs or POD-Galerkin ROMs -- in order to provide time extrapolation.
    
\end{itemize}

\section{$\mu$-POD-LSTM-ROM}
\label{sec:mu_PODLSTMROM}
The first component of the $\mu t$-POD-LSTM-ROM framework is $\mu$-POD-LSTM-ROM, originating from the application of a LSTM autoencoder structure in the context of POD-DL-ROMs. 

While s POD-DL-ROM aims at reducing the dimensionality of the solution by means of a nonlinear projection onto a suitable subspace, the proposed $\mu$-POD-LSTM-ROM framework focuses on the compression of the information necessary to build an entire sequence of solutions.
In particular, a LSTM autoencoder \cite{videolstm} takes a set of sequential inputs and through a LSTM architecture, the \textit{encoder}, provides a lower dimensional representation of the entire sequence as a single vector. Another LSTM based ANN, the \textit{decoder}, takes as input the aforementioned compressed representation and reconstructs the sequence of solutions used to produce it. Hence, such an autoencoder provides a convenient way of representing a sequence of solutions by compressing it in a much lower dimensional single vector, that can be inferred by a properly trained (ANN-based) regressor.

\begin{figure}[h!]
    \centering
    \includegraphics[width=0.9\textwidth]{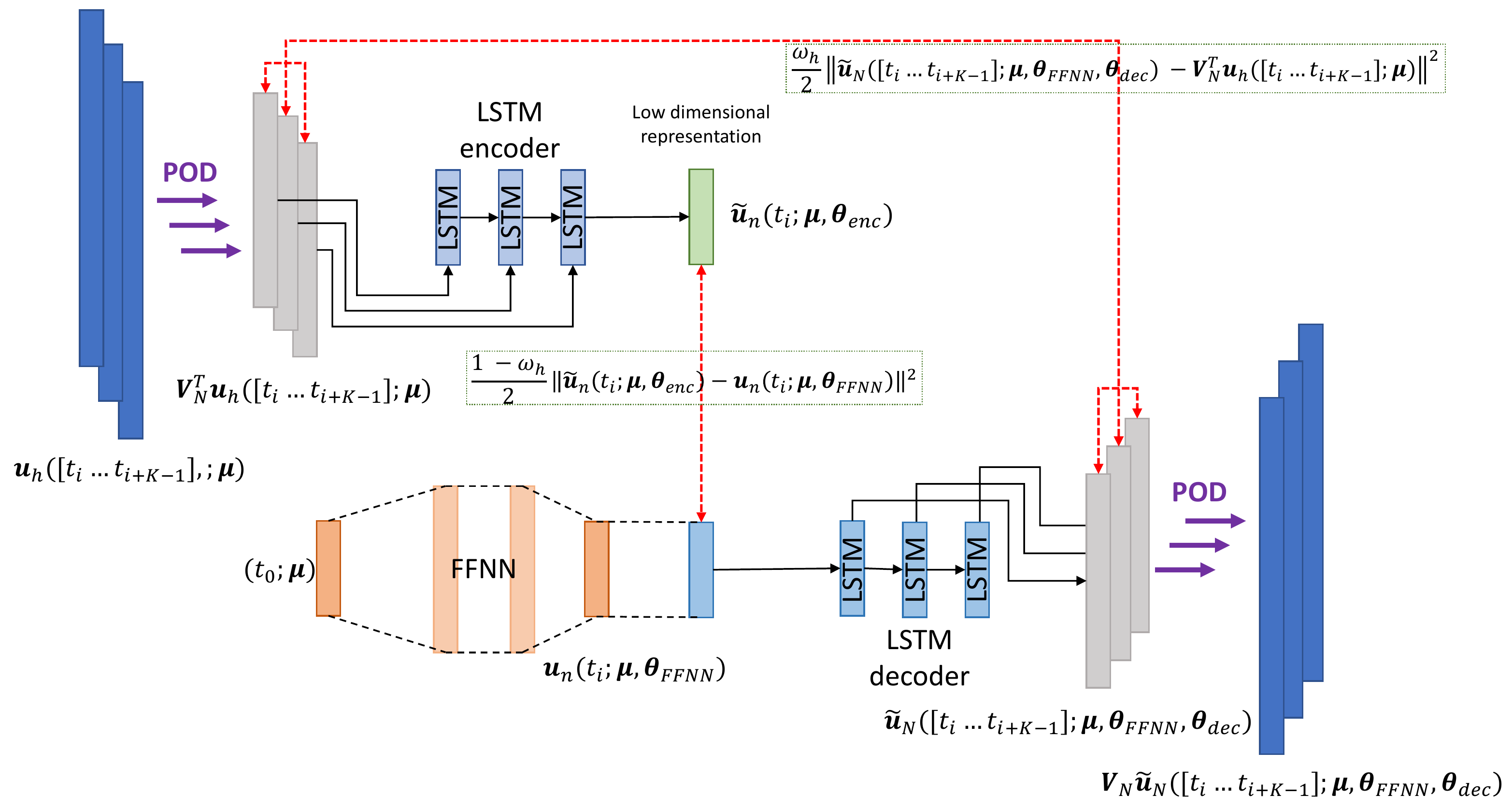}
    \vspace{-3mm}
    \caption{The $\mu$-POD-LSTM-ROM architecture. The full order vectors are reduced by means of POD and assembled in sequences, that are fed to a LSTM autoencoder structure in order to obtain a low-dimensional representation $\mathbf{\tilde{u}}_n(t_i, \boldsymbol{\mu}, \boldsymbol{\theta}_{enc})$ of the entire sequence. A (deep) feedforward neural network is used to infer that low-dimensional representation starting from the initial time $t_i$ and the parameters vector $\boldsymbol{\mu}$ and a LSTM decoder structure allows for the prediction of the sequence of ROM solutions.}
     \vspace{3mm}
    \label{podlstmromfig}
\end{figure}

The logical structure of the $\mu$-POD-LSTM-ROM described above starts with a first dimensionality reduction through the projection of the snapshots onto the POD basis as in a POD-DL-ROM. The snapshots are then put sequentially together in batches and -- optionally -- further reduced via some dense layers. Then, they are passed as sequences to the LSTM encoder architecture and the reconstruction process follows symmetrically through the decoder. At the same time, a feedforward neural network is set up for inferring the hidden representation of the LSTM autoencoder and, from that, for reconstructing the solution sequence starting from the tuple $(\boldsymbol{\mu}, t_i)$\footnote{Note that $t_i \in (0, T - K \Delta t)$ denotes a generic starting time -- not necessarily the one of the solution approximation -- and that the notation $[t_i \dots t_{i+K-1}]$ (with $t_{i+j} = t_i + j \Delta t)$ when used in place of a specific time $t$ indicates the stacking of $K$ subsequent vectors referring to the reported times.}. A schematic representation of this architecture can be found in Figure \ref{podlstmromfig}.

The working scheme of the $\mu$-POD-LSTM-ROM method can be divided in the following blocks:

\begin{itemize}
     \item POD is performed on the snapshot matrix $\boldsymbol{S} \in \mathbb{R}^{N_h \times N_{train} N_t}$ (being $N_h$ the finite discretization dimension of the numerical solver used to produce the snapshots, $N_{train}$ the number of input parameters instances used for building the training set and $N_t$ the number of temporal step used for the time discretization). In the test cases presented, POD was executed in its randomized fashion \cite{halko2011finding}, as in many high-dimensional cases the costs of computing the SVD of the snapshot matrix could become unfeasible \cite{rpodwhy}. This determines a first dimensionality reduction aimed at making the snapshots dimension suitable for feeding the subsequent neural network part. The projection is then performed as 
     \begin{equation}
    \label{proporthdeclstm}
    \mathbf{u}_N(t; \boldsymbol \mu) = \mathbf{V}^T_N \mathbf{u}_h(t; \boldsymbol \mu),
    \end{equation}
    being $\mathbf{V}_N$ the POD projection matrix.

    \item Once dimensionally reduced by means of (r)POD, snapshots are sequentially stacked in matrices of the form $\mathbf{V}^T_N \mathbf{u}_h([t_i \dots t_{i+K-1}]; \boldsymbol \mu) \in \mathbb{R}^{K \times N}$. Such matrices are then grouped in (mini) batches tensors of dimension $\textnormal{dim}_{batch} \times K \times N$ to enable the training.

    \item The sequences of size $K \times N$ from the (mini) batches tensor are then fed to a LSTM encoder structure. The encoder then takes sequentially the $\mathbf{V}^T_N \mathbf{u}_h(t; \boldsymbol \mu)$ as input and modifies its internal state coherently with the evolution of the vectors it receives. The output of the encoder is neglected as it would be of no utility in this context.  
    In particular, the low-dimensional hidden state representation of the LSTM encoder state is built according to the following function:
    \begin{equation}
    \label{hiddenstateeq}
    \Tilde{\mathbf{u}}_n(t_i; \boldsymbol \mu, \boldsymbol{\theta}_{enc}) = {\boldsymbol{\lambda}}_n^{enc}(\mathbf{u}_N([t_i \dots t_{i+K-1}]; \boldsymbol{\mu}, \boldsymbol{\theta}_{enc})).
    \end{equation}
    In the end, the information coming from a sequence of inputs is reduced into a lower dimensional manifold of dimension $n < N \ll N_h$.
    
    \item The reconstruction of the low-dimensional hidden state representation state of the LSTM encoder is performed by a suitable (deep) feedforward neural network consisting in multiple layers of linear transformations and subsequent nonlinear activation functions. The relation learned by this network is 
    \begin{equation}
    \label{reconstructioneq}
    {\mathbf{u}}_n(t_i; \boldsymbol \mu, \boldsymbol{\theta}_{FFNN}) = {\boldsymbol{\phi}}_n^{FFNN}(t_i; \boldsymbol \mu, \boldsymbol{\theta}_{FFNN}).
    \end{equation}
    Note that just the initial time $t_i$ is passed as input to this feedforward neural network, but the hidden representation it infers contains the information coming from an entire sequence of $K$ solutions. This low-dimensional representation is then crucial for the compression of information that allows this architecture to work with sequences, as it provides a convenient way to infer the evolution of the solution by considering a single vector as regression target.
    
    \item The reduced nonlinear trial manifold $\tilde{\mathcal{S}}_N^n$ is modeled using a LSTM decoder that takes as input the approximated hidden representation coming from the feedforward neural network. In particular, the reduced nonlinear trial manifold can be defined as
    \begin{equation}
    \label{nonmanifeq}
    \begin{split}
            \tilde{\mathcal{S}}_N^n = & \{ {\boldsymbol{\lambda}}^{dec}_N(\mathbf{u}_n(t; \boldsymbol{\mu}, {\boldsymbol{\theta}_{FFNN}}); \boldsymbol{\theta}_{dec})_{1 \div N, 1}
    \; | \\ & \; \mathbf{u}_n(t; \boldsymbol{\mu}, \boldsymbol{\theta}_{FFNN}) \in {\mathbb{R}}^{n},    \ t \in [0, T)  \; \textnormal{and} \; \boldsymbol{\mu} \in \mathcal{P} \subset {\mathbb{R}}^{n_{\mu}} \} \subset \mathbb{R}^{N},
    \end{split}
    \end{equation}
    
    and it can be obtained through $\mu t$-POD-LSTM-ROM. Nevertheless, the novel LSTM cell implementation is able to provide also a more meaningful approximated solution manifold to this scope, $\Tilde{\mathcal{S}}_{N,K}^n$, that allows for each input tuple $(\boldsymbol{\mu}, t)$ the reconstruction of the entire sequence $\Tilde{\mathbf{u}}_N([t, \dots, t + (K-1) \Delta t]; \boldsymbol{\mu}) = \mathbf{u}_N([t, \dots, t + (K-1) \Delta t]; \boldsymbol{\mu})$:
     \begin{equation}
    \label{advnonmanifeq}
    \begin{split}
           \tilde{\mathcal{S}}_{N,K}^{n} = & \{ {\boldsymbol{\lambda}}^{dec}_N(\mathbf{u}_n(t; \boldsymbol{\mu}, {\boldsymbol{\theta}_{FFNN}}); \boldsymbol{\theta}_{dec})
    \; | \\ & \; \mathbf{u}_n(t; \boldsymbol{\mu}, \boldsymbol{\theta}_{FFNN}) \in {\mathbb{R}}^{n},    \ t \in [0, T - K \Delta t)  \;  \textnormal{and} \; \boldsymbol{\mu} \in \mathcal{P} \subset {\mathbb{R}}^{n_{\mu}} \} \subset \mathbb{R}^{N \times K}.
    \end{split}
    \end{equation}
    In this context, $\boldsymbol{\lambda}^{dec}_N(\cdot; \boldsymbol{\theta}_{dec}) : {\mathbb{R}}^n \rightarrow {\mathbb{R}^{N \times K}}$ is a suitable LSTM decoder function, taking as input the hidden state of a LSTM encoder and reconstructing the solution starting from it.
     
    \item Once an output sequence has been produced by the decoder function $\boldsymbol{\lambda}^{dec}_N(\cdot; \boldsymbol{\theta}_{dec})$,
    \begin{equation}
    \Tilde{\mathbf{u}}_N([t_i \dots t_{i+K-1}]; \boldsymbol{\mu}, {\boldsymbol{\theta}_{FFNN}}, {\boldsymbol{\theta}_{dec}}) = \boldsymbol{\lambda}^{dec}_N(\mathbf{u}_n(t_i; \boldsymbol{\mu}, {\boldsymbol{\theta}_{FFNN}}); \boldsymbol{\theta}_{dec}),
    \end{equation}
     each of its components is expanded from dimension $N$ to dimension $N_h$ by means of the POD basis found before:
    \begin{equation}
    \Tilde{\mathbf{u}}_h([t_i \dots t_{i+K-1}]; \boldsymbol{\mu}, {\boldsymbol{\theta}_{FFNN}}, {\boldsymbol{\theta}_{dec}}) = \mathbf{V}_N \Tilde{\mathbf{u}}_N([t_i \dots t_{i+K-1}]; \boldsymbol{\mu}, {\boldsymbol{\theta}_{FFNN}}, {\boldsymbol{\theta}_{dec}});
    \end{equation}
    finally, a stack of $N_h$ full dimensional time sequential solutions is found.

 \end{itemize}
 
  In the notation used above, parameters vectors $\boldsymbol{\theta}_{enc}$, $\boldsymbol{\theta}_{FFNN}$ and $\boldsymbol{\theta}_{dec}$ contain the trainable parameters of the networks. Their hyperparameters (such as, \textit{e.g.}, the number of stacked LSTM cells or their possible bidirectionality or the depth of the DFNN) should be considered as well in a separate optimization process.
  
\begin{remark}
Note that, optionally, the dimensionality of each POD-reduced snapshot composing the input sequence can be further reduced by means of a time distributed feedforward neural network before passing through the LSTM autoencoder. This network applies to each $\mathbf{u}_N(t; \boldsymbol \mu) = \mathbf{V}^T_N \mathbf{u}_h(t; \boldsymbol \mu)$ in the sequence $ \mathbf{u}_N([t_i \dots t_{i+K-1}]; \boldsymbol \mu)$ the same nonlinear transformation that reduces its dimensions from $N$ to $N_{red}$. For simplicity, in this work we will assume $N_{red} = N$. Furthermore, also the low-dimensional representation provided by the autoencoder can be optionally further compressed by means of a feedforward neural network.
 \end{remark}

\begin{remark}
 Note that $N_h$ and $N$ represent just discretization dimensions. In case of vectorial PDE problems, for each time at each spatial discretization point we associate a vector in $\mathbb{R}^{n_{ch}}$. In this case the dimensionality of the FOM solution increases from $N_h$ to $N_h n_{ch}$ and the one of the POD reduced solution increases from $N$ to $N n_{ch}$. The structure just described, though, is still valid also in this case provided that the involved dimensions are suitably modified.
 \end{remark}
 
 The offline training stage consists of the solution of an optimization problem in which a loss function expressed as a function in the variable $\boldsymbol{\theta} = (\boldsymbol{\theta}_{enc}, \boldsymbol{\theta}_{FFNN}, \boldsymbol{\theta}_{dec})$ should be minimized.  
 In particular, for the training of the $\mu$-POD-LSTM-ROM, the snapshot matrix $\mathbf{S} \in \mathbb{R}^{N_h \times N_{train} N_t}$ (with $N_{train}$ being the number of unique instances drawn from the parameters' space and $N_t$ is the number of timesteps chosen for the time discretization of the interval (0, T)) is compressed by means of POD as explained before to become $\mathbf{S}_{POD} \in \mathbb{R}^{N \times N_{train} N_t}$. Then, sequences from this matrix are extracted to form the so called \textit{base tensor} $\mathbf{T} \in \mathbb{R}^{N_{train}(N_t - K) \times N \times K}$, to be fed to the network,
 \begin{equation*}
     \mathbf{T}(i, j, k) = (\mathbf{u}_N(t_{\alpha_i} + k \Delta t, \boldsymbol{\mu}_{\beta_i}))_j,
 \end{equation*}
 with $\alpha_i = i \textnormal{mod} (N_t - K)$ and $\beta_i = \frac{i - \alpha_i}{N_t - K}$ and $(\cdot)_j$ denoting the extraction of the $j^{th}$ component from a vector. 
 
The minimization problem can therefore be formulated in this case as 
 \begin{equation}
\min_{\boldsymbol{\theta}} \mathcal{J}(\boldsymbol{\theta}) = \min_{\boldsymbol{\theta}} \frac{1}{N_{train}(N_t - K)}\sum_{i=1}^{N_{train}} \sum_{k=1}^{N_t - K} \mathcal{L}(t_k, \boldsymbol \mu_i; \boldsymbol{\theta}),
\label{minimization_problem_LSTM}
\end{equation}
where the loss function is defined by
\begin{equation}
\mathcal{L}(t_k, \boldsymbol{\mu}_i;  {\boldsymbol{\theta}}) = \frac{\omega_h}{2} \mathcal{L}_{rec}(t_k, \boldsymbol{\mu}_i;  {\boldsymbol{\theta}}) + \frac{1-\omega_h}{2}       \mathcal{L}_{hid}(t_k, \boldsymbol{\mu}_i;  {\boldsymbol{\theta}}),
\label{eq:loss_encoder}
\end{equation}
with
\[
\mathcal{L}_{rec}(t_k, \boldsymbol{\mu}_i;  {\boldsymbol{\theta}}) = \| \mathbf{T}(i,:,k) - \mathbf{\tilde{u}}_N(t_{\alpha_i} + k \Delta t; \boldsymbol{\mu}_{\beta_i},  {\boldsymbol{\theta}_{FFNN}, \boldsymbol{\theta}_{dec}})\|^2
\]
and
\[
\mathcal{L}_{int}(t^k, \boldsymbol{\mu}_i;  {\boldsymbol{\theta}}) =  \| \tilde{\mathbf{u}}_n(t_{\alpha_i} + k \Delta t; \boldsymbol{\mu}_{\beta_i}, \boldsymbol{\theta}_{enc}) -  {\mathbf{u}}_n(t_{\alpha_i} + k \Delta t; \boldsymbol{\mu}_{\beta_i},  {\boldsymbol{\theta}_{FFNN}})\|^2.
\]
The loss function (\ref{eq:loss_encoder}) penalizes the reconstruction error from the LSTM autoencoder through $\mathcal{L}_{rec}$ and the difference between the low-dimensional hidden representation learned by the encoder and the prediction from the feedforward neural network fed with the problem parameters through $\mathcal{L}_{int}$. The coefficient $\omega_h \in [0,1]$ regulates the importance of the two components of the loss function.

During the online stage (at testing), just the feedforward part $\boldsymbol{\phi}^{FFNN}_N(\cdot; \cdot, \boldsymbol{\theta}_{FFNN})$ and the decoder $\boldsymbol{\lambda}_N^{dec}(\cdot; \boldsymbol{\theta}_{dec})$ are used. The encoder part is added at training time in order to help the network learning the correct hidden representation of the sequences in a data-driven fashion.

\section{$t$-POD-LSTM-ROM}
\label{sec:t_PODLSTMROM}
$t$-POD-LSTM-ROM is the second component of the $\mu t$-POD-LSTM-ROM framework, providing it with time extrapolation capabilities. As the name suggests, it works on time series forecasting by solving iteratively a \textit{sequence to sequence} problem.

We define a \textit{sequence to sequence} (also referred to as \textit{seq2seq}) learning problem (see, e.g., \cite{sangiorgio2020}) as the forecasting of a certain number ($k$) of steps ahead in a time series $y(t)$. Therefore, the solution of the problem is a model (the \textit{predictor}), that takes as input $p$ time steps in the past and returns as output the forecasted $k$ steps ahead in the future. As a matter of fact, this latter is a function taking as input a sequence and allowing to predict an output one, which can be summarized by the following expression:    
\begin{equation}
    \Phi(\cdot) : [y(t-p+1), y(t-p+2), \dots, y(t)] \longrightarrow [\hat{y}(t+1), \hat{y}(t+2), \dots, \hat{y}(t+k)].    
\end{equation}

Traditional machine learning models, such as simple regression, support vector regression, ARIMA and feedforward neural networks have been used to tackle the problem \cite{pastimeseriesDu}. Hidden Markov models or fuzzy logic based models have also proven to be somehow effective in the field \cite{hmm4ts,fuzzy4ts}. 
Recently though, artificial neural networks featuring recurrence mechanisms such as simple recurrent neural networks or LSTM cells have become the standard for time series prediction when having large amount of data available for the training \cite{pastimeseriesDu,schmid2001,sangiorgio2020,arima4ts}.

The idea of a windowed auto-regressive prediction has been exploited in neural ODEs \cite{chen2019neural,massaroli2021dissecting}, where deterministic numerical solvers consider also statistically learned residuals in order to perform the PDE integration in time. This method is very effective in time extrapolation capabilities, but still considers numerical integration of high dimensional systems. A similar approach exploited in the field of ROMs has been considered in \cite{osti20}, where a correction parameter was included in the numerical time integration process. Since the use of a numerical integrator can negatively impact on time performances of the method, we decided to rely on a time series forecasting problem as described before. In particular, in \cite{lstmforpdes} a time series approach for PDE problems using LSTMs proved to be effective in forecasting reduced barotropic climate models and showing interesting time performances, however  neglecting in that case  the parametric nature of the problem. Note that, also for this case, considering architectures based on the \textit{seq2seq} paradigm shows close similarities with the behavior of numerical solvers, as they build predictions for future times based on the past, mimicking time marching numerical schemes. \\

In this work, we introduce the $t$-POD-LSTM-ROM architecture, aimed at solving the \textit{seq2seq} problem in the context of ROMs for parameterized dynamical systems. In particular, the novel framework allows to extend in the temporal dimension the ROM solution provided by the $\mu$-POD-LSTM-ROM framework introduced before, by considering the \textit{seq2seq} problem on the reduced order vectors. The training of this additional neural network does not require an increased number of snapshots, keeping the temporal cost of the offline phase relatively low. Furthermore, the possibility to extend the temporal domain of definition of the solution could in principle allow for the training using full order snapshots defined on a shorter period of time, thus requiring less computational expenses when generating them.

\begin{figure}[t]
  \vspace{-0.2cm}
    \centering
    \includegraphics[width=0.9\textwidth]{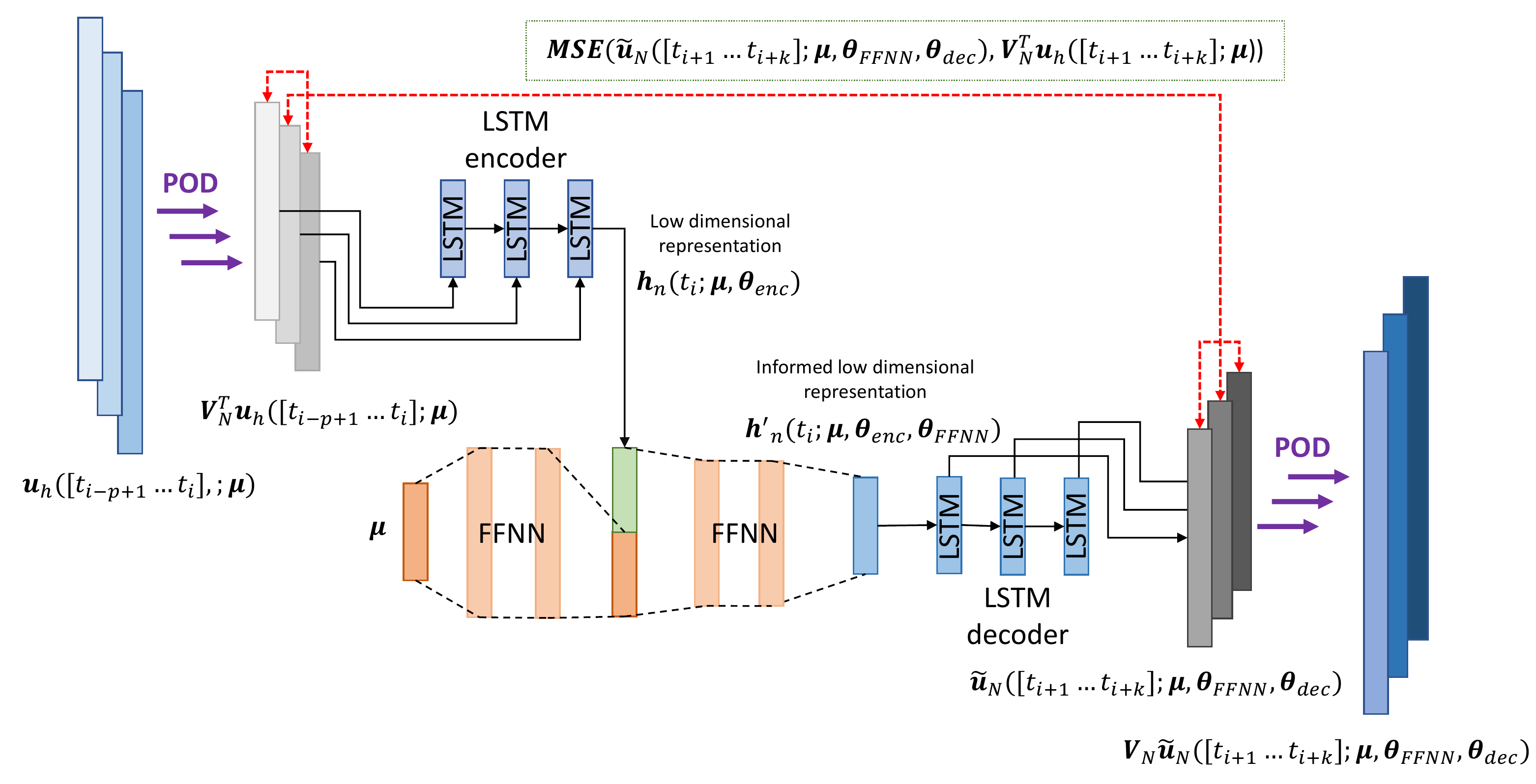}
    \vspace{-0.2cm}
    \caption{The $t$-POD-LSTM-ROM architecture. After an initial POD-based dimensionality reduction, the sequence of reduced vectors passes through a LSTM encoder. The low-dimensional representation $\mathbf{h}_n(t_i; \boldsymbol{\mu}, \boldsymbol{\theta}_{enc})$ thus obtained is then enriched with the information coming from the problem's parameters and starting from it a LSTM decoder finally extracts the $k$-steps forecast sequence.}
    \label{lstmtsfig}
\end{figure}

The architecture of the $t$-POD-LSTM-ROM, summarized in Figure \ref{lstmtsfig}, consists of the components listed below:

\begin{itemize}
    \item The same POD-reduced sequences $\mathbf{V}_N^T\mathbf{u}_h([t_{i-p+1} \dots t_i]; \boldsymbol{\mu})$ created to feed the $\mu$-POD-LSTM-ROM framework, and formed by the $p-1$ snapshots preceding the one at time $t_i$ and the one at time $t_i$ itself, are collected in pairs with $\mathbf{V}_N^T\mathbf{u}_h([t_{i+1} \dots t_{i+k}]; \boldsymbol{\mu})$, the sequence formed by the $k$ following reduced snapshots to be used in the training phase.
    
    \item A LSTM autoencoder is then applied to provide a low-dimensional representation of the data extracted by the past time steps sequence. In particular, the encoder acts according to
    \begin{equation}
        \mathbf{h}_n(t_i; \boldsymbol{\mu}, \boldsymbol{\theta}_{enc}) = \boldsymbol{\Lambda}_n^{enc}(\mathbf{V}_N^T\mathbf{u}_h([t_{i-p+1} \dots t_i]; \boldsymbol{\mu}); \boldsymbol{\theta}_{enc}).
    \end{equation}
    The low-dimensional representation then contains in principle all the necessary information to build the sequence of future solutions by means of a LSTM decoder.
    
    \item A feedforward neural network is then used to inform the architecture on the parameters of the system. First, some dense layers are applied in order to expand the information contained in the parameters' vector $\boldsymbol{\mu}$ and then the result of this operation is concatenated to the low-dimensional representation as $\mathbf{h}_n(t_i; \boldsymbol{\mu}, \boldsymbol{\theta}_{enc}) \oplus \boldsymbol{\phi}(\boldsymbol{\mu}, \boldsymbol{\theta}_{FFNN1})$\footnote{We define in this context $\oplus$ with the concatenation between two vectors by appending the second after the first. By defining $\mathbf{a} \in \mathbb{R}^{n_1}$ and $\mathbf{b} \in \mathbb{R}^{n_2}$, then $\mathbf{a} \oplus \mathbf{b} = [\mathbf{a}, \mathbf{b}] = \mathbf{c} \in \mathbb{R}^{n_1+n_2}$}. This concatenation is then fed to another feedforward neural network aimed at merging the information coming from the past snapshots and the one coming from the parameters in order to form a better low-dimensional representation of the information gathered so far, of the form
    \begin{equation}
        {\mathbf{h'}}_n(t_i; \boldsymbol{\mu}, \boldsymbol{\theta}_{enc}, \boldsymbol{\theta}_{FFNN}) = {\boldsymbol{\phi'}}(\boldsymbol{\mu}, [\mathbf{h}_n(t_i; \boldsymbol{\mu}, \boldsymbol{\theta}_{enc}) \oplus \boldsymbol{\phi}(\boldsymbol{\mu}, \boldsymbol{\theta}_{FFNN1})]; \boldsymbol{\theta}_{FFNN2}),
    \end{equation}
    with $\boldsymbol{\theta}_{FFNN} = (\boldsymbol{\theta}_{FFNN1}, \boldsymbol{\theta}_{FFNN2})$, being $\boldsymbol{\theta}_{FFNN1}$ and $\boldsymbol{\theta}_{FFNN2}$ the vectors of parameters of the two feedforward neural network parts, namely $\boldsymbol{\phi}$ and $\boldsymbol{\phi'}$.
    
    \item A LSTM decoder then takes as input the informed low-dimensional representation, given by ${\mathbf{h'}}_n(t_i; \boldsymbol{\mu}, \boldsymbol{\theta}_{enc}, \boldsymbol{\theta}_{FFNN})$, and extracts the forecast reduced future sequence according to
    \begin{equation}
        \mathbf{\Tilde{u}}_N([t_{i+1} \dots t_{i+k}]; \boldsymbol{\mu}, \boldsymbol{\theta}_{FFNN}, \boldsymbol{\theta}_{enc}, \boldsymbol{\theta}_{dec}) = \boldsymbol{\Lambda}_n^{dec}({\mathbf{h'}}_n(t_i; \boldsymbol{\mu}, \boldsymbol{\theta}_{enc}, \boldsymbol{\theta}_{FFNN})).
    \end{equation}
    From this output, the sequence of full order predicted solutions is then reconstructed by means of POD basis. 
    In contrast with the $\mu$-POD-LSTM-ROM technique considered before (and the POD-DL-ROM one), in this case the architecture of the autoencoder -- and therefore that of the entire network -- remains the same during both the training and the testing stages.
\end{itemize}

Also in this case, the parameters vectors $\boldsymbol{\theta}_{enc}$, $\boldsymbol{\theta}_{FFNN}$ and $\boldsymbol{\theta}_{dec}$ group the trainable weights and biases of the networks. 
The training of the network described above is then performed as an optimization problem in which a loss function expressed in the variable $\boldsymbol{\theta} = (\boldsymbol{\theta}_{enc}, \boldsymbol{\theta}_{FFNN}, \boldsymbol{\theta}_{dec})$ should be minimized. After the definition of the \textit{base tensor} $\mathbf{T} \in \mathbb{R}^{N_{train}(N_t - K) \times N \times K}$, as described when introducing $\mu$-POD-LSTM-ROM, with $K = k + p$, we divide it in two parts: the \textit{previous steps} tensor $\mathbf{P} \in \mathbb{R}^{N_{train}(N_t - K) \times N \times p}$ and the \textit{horizon} tensor $\mathbf{H} \in \mathbb{R}^{N_{train}(N_t - K) \times N \times k}$. In particular, we define 
\begin{equation}
    \mathbf{P} = \mathbf{T}(:, :, 1:p) \textnormal{ and } \mathbf{H} = \mathbf{T}(:, :, (p+1):K),
\end{equation}
so that the tensor $\mathbf{P}$ represents the previous time steps to be used for the prediction and $\mathbf{H}$ contains the target sequences to forecast.

In the end, the minimization problem solved during training can be defined as 

 \begin{equation}
\min_{\boldsymbol{\theta}} \mathcal{J}(\boldsymbol{\theta}) = \min_{\boldsymbol{\theta}} \frac{1}{N_{train}(N_t - K)}\sum_{i=1}^{N_{train}} \sum_{j=1}^{N_t - K} \mathcal{L}(t_j, \boldsymbol \mu_i; \boldsymbol{\theta}),
\label{minimization_problem_LSTM_1}
\end{equation}

\noindent
where we define \begin{equation}
\mathcal{L}(t_j, \boldsymbol{\mu}_i;  {\boldsymbol{\theta}}) = \textnormal{MSE}\left[\Tilde{\mathbf{u}}_N([t_{j+1} \dots t_{j+k}]; \boldsymbol{\mu}_i, \boldsymbol{\theta}_{FFNN}, \boldsymbol{\theta}_{enc}, \boldsymbol{\theta}_{dec}),\ \mathbf{V}_N^T\mathbf{u}_h([t_{j+1} \dots t_{j+k}]; \boldsymbol{\mu}_i)\right]
\label{eq:loss_lstmts}
\end{equation}
with
\begin{equation*}
\begin{split}
    \textnormal{MSE}&\left[\Tilde{\mathbf{u}}_N([t_{j+1} \dots t_{j+k}]; \boldsymbol{\mu}, \boldsymbol{\theta}_{FFNN}, \boldsymbol{\theta}_{enc}, \boldsymbol{\theta}_{dec}),\ \mathbf{V}_N^T\mathbf{u}_h([t_{j+1} \dots t_{j+k}]; \boldsymbol{\mu})\right] = \\ & \qquad \qquad \frac{1}{N k}\sum_{l=1}^{N}\sum_{p=1}^{k}\left(\left(\Tilde{\mathbf{u}}_N(t_{j+p}; \boldsymbol{\mu}_i, \boldsymbol{\theta}_{FFNN}, \boldsymbol{\theta}_{enc}, \boldsymbol{\theta}_{dec})\right)_l - \left(\mathbf{V}_N^T\mathbf{u}_h(t_{j+p}; \boldsymbol{\mu}_i)\right)_l\right)^2.
    \end{split}
\end{equation*}

\noindent The loss function (\ref{eq:loss_lstmts}) therefore penalizes prediction errors and maximizes the accuracy in prediction. 

\section{Results}
\label{sec:results}

In this section we present a set of numerical results obtained on three different test cases, related with {\em (i)} a 3 species Lotka-Volterra equations (Section \ref{LVsection}), {\em (ii)} unsteady advection-diffusion-reaction equation (Section \ref{ADRsection}), {\em (iii)} incompressible Navier-Stokes equations (Section \ref{NSsection}). To assess the accuracy of the numerical results, we consider the same two error indicators defined in \cite{fresca2021poddlrom}, namely:
\begin{itemize}
    \item the error indicator $\epsilon_{rel} \in \mathbb{R}$ defined as
\begin{equation}
\epsilon_{rel}(\mathbf{u}_h, \mathbf{\tilde{u}}_h) = \frac{1}{N_{test}} \sum_{i  = 1}^{N_{test}} \left(\displaystyle \frac{\sqrt{ \sum_{k=1}^{N_t} || \mathbf{u}^k_h(\boldsymbol{\mu}_{test,i}) - \mathbf{\tilde{u}}^k_h(\boldsymbol{\mu}_{test,i}) ||^2}}{\sqrt{\sum_{k=1}^{N_t} || \mathbf{u}_h^k(\boldsymbol{\mu}_{test,i}) ||^2}} \right),
\label{eq:error_indicator}
\end{equation}
\item the relative error $\boldsymbol{\epsilon}_k \in \mathbb{R}^{\sum_i N_h^i}$, for $k = 1, \ldots, N_t$, defined as
\begin{equation}
\displaystyle \boldsymbol{\epsilon}_k(\mathbf{u}_h, \mathbf{\tilde{u}}_h) = \displaystyle \frac{ | \mathbf{u}^k_h(\boldsymbol{\mu}_{test}) - \mathbf{\tilde{u}}^k_h(\boldsymbol{\mu}_{test}) |}{\sqrt{\frac{1}{N_t}\sum_{k=1}^{N_t} || \mathbf{u}^k_h(\boldsymbol{\mu}_{test}) ||^2}}.
\label{eq:relative_error}
\end{equation}
\end{itemize}
Note that the error indicator $\epsilon_{rel}$ provides a (scalar) numerical estimation of the accuracy performances of the method on the entire test set.

The $\mu$-POD-LSTM-ROM and the $\mu t$-POD-LSTM-ROM architectures have been developed using TensorFlow 2.4 framework \cite{tensorflow2015-whitepaper}. FOM data for test cases (ii) and (iii) have been obtained by the redbKIT v2.2 library \cite{reporedbkit}, implementing the methods described in \cite{quarteroni2016reduced}.  All the simulations have been run on an Intel\textregistered \ Core i9 @ 2.40GHz CPU, 16 GB RAM and NVIDIA\textregistered \ GTX1650 video card personal computer. 

\subsection{Lotka-Volterra competition model (3 species)}
\label{LVsection}
The first test case is the 3 species Lotka-Volterra competition model, selected to provide a proof-of-concept of the method on a simple but aperiodic test case.
The goal is the reconstruction of the solution $\mathbf{u} = \mathbf{u}(t; \mu) \in \mathbb{R}^3$ of the following system:
\begin{equation}
    \left\{
\begin{aligned}
& \frac{d u_1}{dt}(t) = u_1(t)(\mu - 0.1u_1(t) - 0.5u_2(t) - 0.5u_3(t)) & \  & t \in (0,T),\\
& \frac{d u_2}{dt}(t) = u_2(t)(-\mu + 0.5u_1(t) - 0.3u_3(t)) & \  & t \in (0,T) \\
& \frac{d u_3}{dt}(t) = u_3(t)(-\mu + 0.2u_1(t) + 0.5u_2(t)) & \  & t \in (0,T) \\
& u_i(0) = 0.5 & \  &  \forall i \in \{1,2,3\}.
\end{aligned}
\right.
\label{eqLV}
\end{equation}

Note that due to the low-dimensionality of the problem ($N_h = 3$), in this case the use of POD is not necessary and therefore it is not performed. Nevertheless, the entire deep learning-based architecture is still used, providing a first glance on the performances of the presented framework when considering time extrapolation capabilities.

The parameter $\mu \in \mathcal{P} = [1, 3]$ models both the reproduction rate of the species 1 (the prey) and the mortality rate of species 2 and 3 (predators), assumed to be equal. The impact of $\mu$ on the solution regards both the amplitude and the frequency of the oscillation of the 3 species' populations. \\
Equations have been discretized by means of an explicit Runge-Kutta (4,5) formula considering a time step $\Delta t = 0.1$ over the time interval $(0, T)$, with $T=9.9$.

The LSTM-ROM framework used to find the solution of the system considers $N_t = 100$ time instances with $N_{train} = 21$. In particular, the selected $\mu$ for the training are equally spaced in the interval $\mathcal{P} = [1, 3]$ (that is, $\mathcal{P}_{train} = \{1, 1.1, 1.2, \dots, 2.9, 3\}$). The LSTM sequence length used for the training is $K=20$, the hidden dimension of the LSTM network has been chosen to be $n = 40$ and the loss parameter $\omega_h$ introduced in (\ref{eq:loss_encoder}) has been set equal to $\omega_h = 0.9$. These choices are the result of a random search hyperparameters tuning \cite{bergstra12} considering both accuracy and time performances. The number of epochs have been fixed to a maximum of $n_{epochs} = 4000$ with the early stopping criterion intervening after 50 epochs of missed improvement of the loss function over the validation set during optimization. \\

We present here, for the sake of comparison, also the results obtained with the DL-ROM framework, for which we considered the same training set used to train the $\mu$-LSTM-ROM network with the same maximum number of epochs and early stopping criterion. The low-dimensional manifold representation was $n = 40$ for $\mu$-LSTM-ROM and $n=10$ for DL-ROM. Note that in this case $n > N_h$. This happens in the context of $\mu$-LSTM-ROM because the low-dimensional representation should contain enough information for the decoder to reconstruct an entire sequence of $N_h = 3-$dimensional vectors of length $K=20$ resulting in a $N_h \times K = 3 \times 20 = 60$ components output. Also DL-ROM in this case showed better performances with $n=10 > n_{\boldsymbol{\mu}}+1$, probably due to the extremely low magnitude of $N_h$. In general though, the proposed LSTM-based framework requires larger dimensional reduced manifolds with respect to POD-DL-ROM in order to provide the decoder with enough information to reconstruct POD-reduced solution sequences. 

\begin{figure}[t]
\begin{center}
    \large \textbf{Simulated results - $\boldsymbol{\mu}$ = 1.95}
\end{center}

\begin{minipage}[]{0.32\textwidth}
    \includegraphics[width=4.1cm]{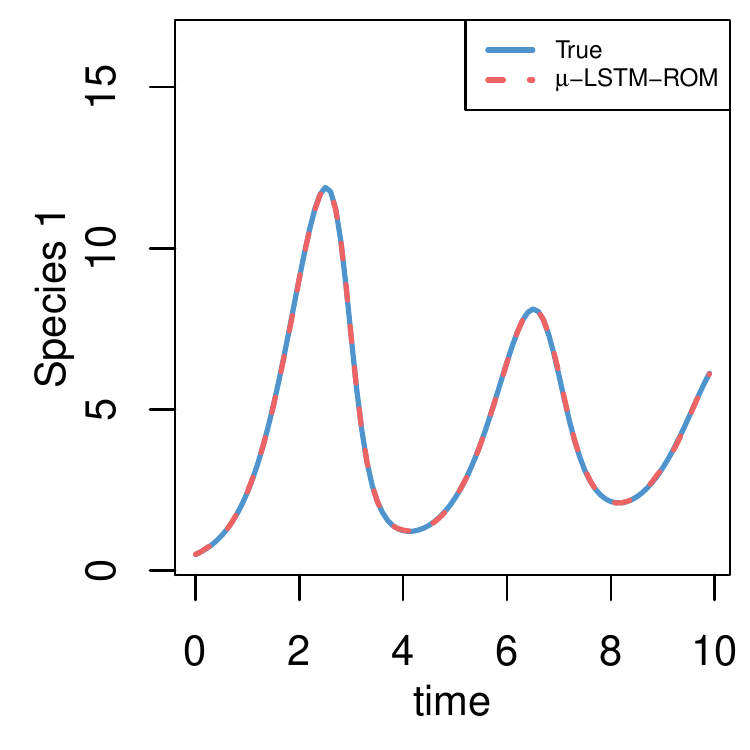}\\
\end{minipage}
\begin{minipage}[]{0.32\textwidth}
    \includegraphics[width=4.1cm]{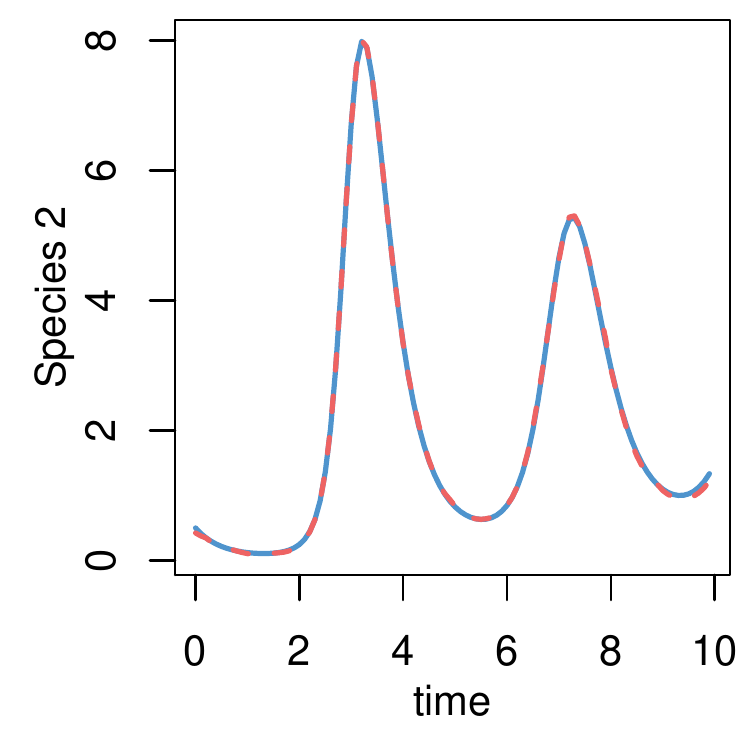}\\
\end{minipage}
\begin{minipage}[]{0.32\textwidth}
    \includegraphics[width=4.1cm]{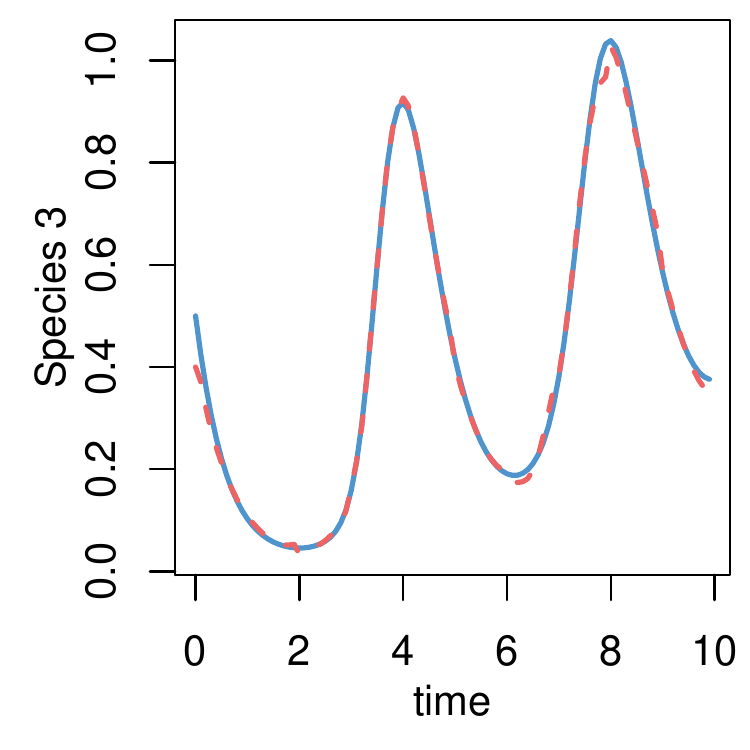}\\
\end{minipage}
\vspace{-0.6cm}

\begin{minipage}[]{0.32\textwidth}
    \includegraphics[width=4.1cm]{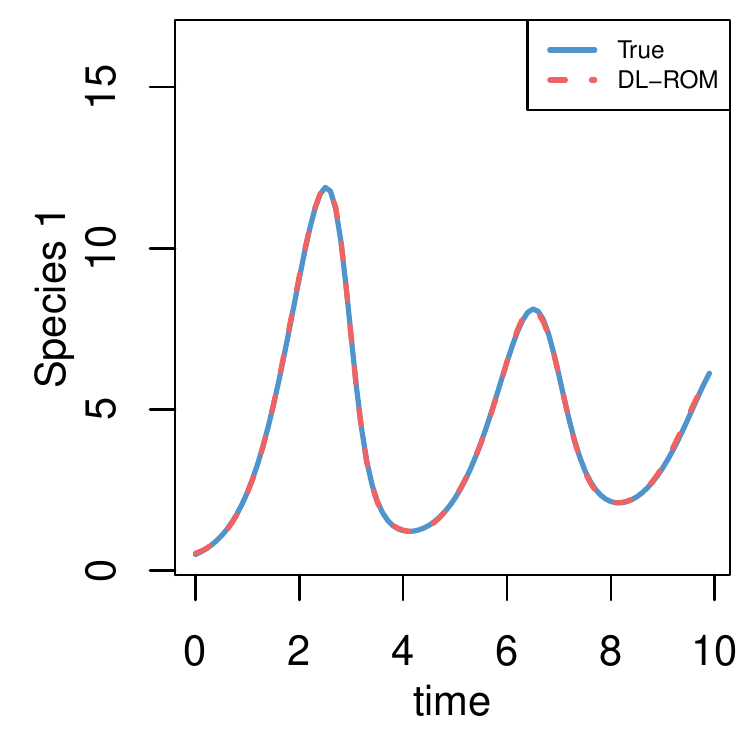}\\
\end{minipage}
\begin{minipage}[]{0.32\textwidth}
    \includegraphics[width=4.1cm]{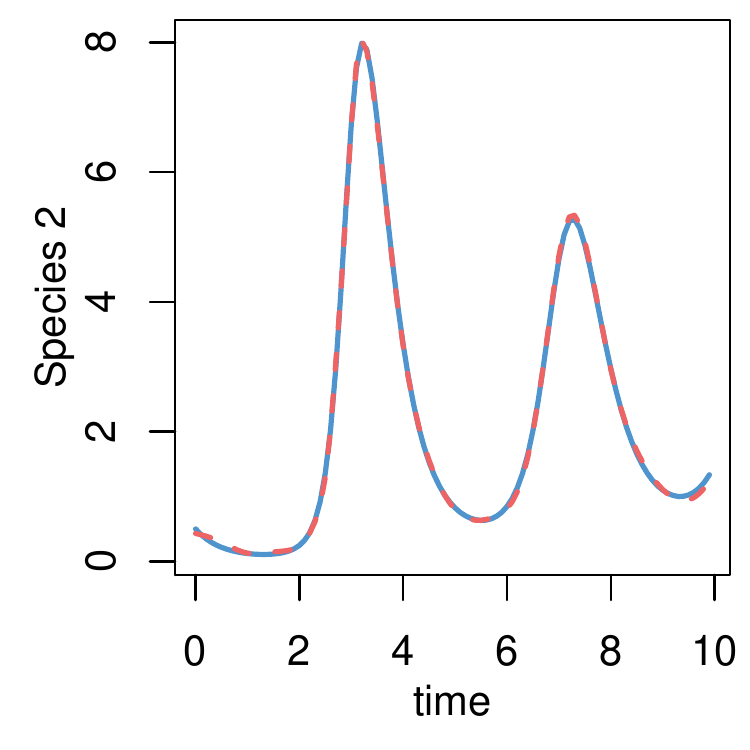}\\
\end{minipage}
\begin{minipage}[]{0.32\textwidth}
    \includegraphics[width=4.1cm]{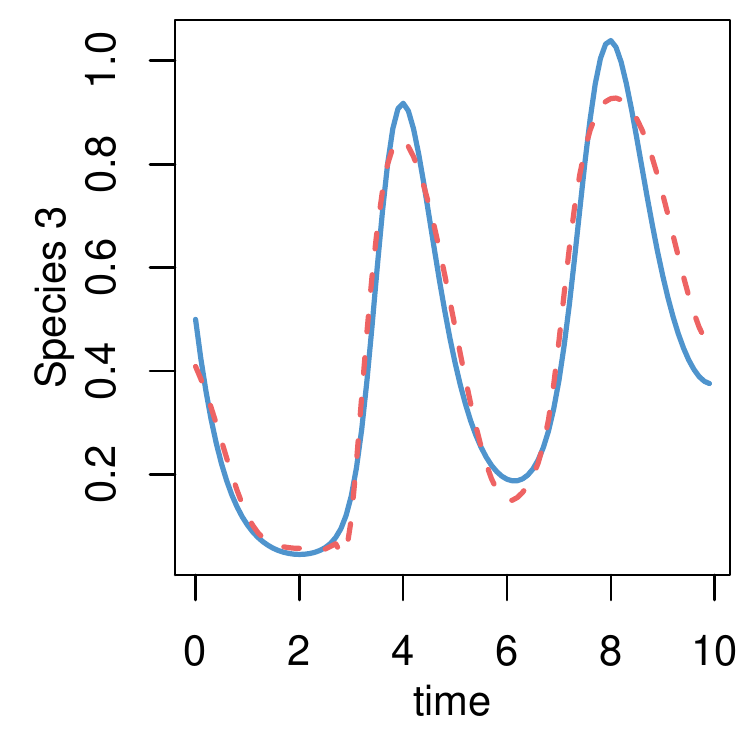}\\
\end{minipage}
\vspace{-0.6cm}
\caption{Test case 1 -- Lotka-Volterra system. Simulated results for $\mu = 1.95$. Top: $\mu$-LSTM-ROM framework, bottom: DL-ROM framework.}
\label{figLV}
\end{figure}

Regarding the architectures considered, we chose to rely on a much larger LSTM-ROM architecture (34433 trainable parameters) with respect to the DL-ROM one used (2943 trainable parameters). This unbalance is recurrent in all the test cases presented. $\mu$-POD-LSTM-ROM requires in fact a larger number of parameters as it needs to encode more information with respect to POD-DL-ROM. $\mu$-LSTM-ROM training took 3124 epochs (932s), while DL-ROM one took 1441 epochs (274s). Results in terms of time evolution of the 3 species for a representative instance of parameters space ($\mu = 1.95$, equally distant from the extremes of $\mathcal{P}_{test}$) are reported in Figure \ref{figLV}.

The error indicator $\epsilon_{rel}^{LSTM-ROM}$ for the LSTM-ROM case is $5.582 \cdot 10^{-3}$, while for the DL-ROM we find $\epsilon_{rel}^{DL-ROM} = 1.363 \cdot 10^{-2}$. Hence the novel framework provides slightly better results in terms of accuracy in this particular application. Table \ref{table1} reports the obtained relative error $\boldsymbol{\epsilon}_k$.

\begin{table}[h]
\centering
\begin{tabular}{@{}lll@{}}
\cmidrule(l){2-3}
                  & $\mathbf{\boldsymbol{\epsilon}_{k}^{mean}}$ & $\mathbf{\boldsymbol{\epsilon}_{k}^{max}}$ \\ \midrule
\textbf{$\mu$-LSTM-ROM} & $4.836 \cdot 10^{-4}$                & $6.703 \cdot 10^{-3}$                \\ \midrule
\textbf{DL-ROM}   & $1.076 \cdot 10^{-3}$                  & $1.082 \cdot 10^{-2}$                 \\ \bottomrule
\end{tabular}
\caption{Test case 1 -- Lotka-Volterra system. Error results in comparison between $\mu$-LSTM-ROM and DL-ROM methods.}
\label{table1}
\end{table}

The $95\%$ bootstrap confidence intervals for the mean relative error $\mathbf{\boldsymbol{\epsilon}_{k}^{mean}}$ in the LSTM-ROM case is given by $C_{0.95}^{LSTM-ROM} = [4.660 \cdot 10^{-4}, 5.007 \cdot 10^{-4}]$, while in the DL-ROM case the same confidence interval is $C_{0.95}^{LSTM-ROM} = [1.040 \cdot 10^{-3}, 1.110 \cdot 10^{-3}]$. In particular, the loss in accuracy regarding both methods when considering the third species has to be accounted to the smaller scales characterizing it. Furthermore, we need to consider the disparity in the number of parameters among the two neural networks. Both architectures though are extremely simple and the DL-ROM one includes the deepest networks and the most pronounced nonlinearities; hence these features should partly compensate for its lower parameters' number.

Regarding prediction times, we consider them from a double perspective. We will report the crude time $t_{NN}$ that the neural network structure takes to perform the forward pass and thus produce a prediction of the solution for each instance of the parameters' test space $\mathcal{P}_{test}$, and the total time $t_{rec}$ which considers also the construction of the full order vectors. Note that the $\mu$-POD-LSTM-ROM neural network would require $\frac{1}{K}$ input-output queries with respect to the POD-DL-ROM one in this phase, because it needs just one prediction every K time steps to be able to build the time evolution of the solution. The quantity $t_{rec}$ takes also into account the amount of time required to assemble the actual solutions, but this is strictly code and language dependent. The bottleneck in computations is in our opinion represented by $t_{NN}$.

\noindent
For the case at stake we found, after a run of 100 queries on the entire $\mathcal{P}_{test}$, the times reported in Table \ref{table2}. A 34.6\% reduction in $t_{NN}^{mean}$ and a 31.0\% reduction in $t_{rec}^{mean}$ can be observed when using $\mu$-LSTM-ROM instead of DL-ROM, thus highlighting an increased efficiency of the new method with respect to the old one also when considering the entire solution reconstruction phase.

\begin{table}[h]
\centering
\begin{tabular}{@{}llll@{}}
\cmidrule(l){2-4}
                  & $\mathbf{t_{NN}^{min}}$ & $\mathbf{t_{NN}^{mean}}$ & $\mathbf{t_{NN}^{max}}$ \\ \midrule
\textbf{$\mu$-LSTM-ROM} & 0.0364s               & 0.0415s                & 0.0598s               \\ \midrule
\textbf{DL-ROM}   & 0.0583s                & 0.0635s                 & 0.0813s                \\ \bottomrule
\end{tabular}
\\

\vspace{0.5cm}

\centering
\begin{tabular}{@{}llll@{}}
\cmidrule(l){2-4}
                  & $\mathbf{t_{rec}^{min}}$ & $\mathbf{t_{rec}^{mean}}$ & $\mathbf{t_{rec}^{max}}$ \\ \midrule
\textbf{$\mu$-LSTM-ROM} & 0.0381s               & 0.0446s                & 0.3281s               \\ \midrule
\textbf{DL-ROM}   & 0.0587s                & 0.0646s                 & 0.0802s                \\ \bottomrule
\end{tabular}
\caption{Test case 1 -- Lotka-Volterra system. Temporal results for the comparison between $\mu$-LSTM-ROM and DL-ROM frameworks.}
\label{table2}
\end{table}

\subsubsection*{Lotka-Volterra time extrapolation}
We then tested the time extrapolation capabilities of the method both in the short and the long term. In particular, we present the results of the $\mu t$-LSTM-ROM framework for the same training and test sets ($\mathcal{P}_{train}$ and $\mathcal{P}_{test}$) just considered. In this case, though, the training snapshots have been acquired from the time interval $(0,T)$, with $T=9.9$ and $\Delta t = 0.1$, while the testing ones are taken in $[T_{in}, T_{fin}]$, where $T_{in} = 5$ and $T_{fin} = 14.9$. Training parameters are the same as in the interpolation case. We therefore consider an extrapolation window of length $T_{ext}=5$, half of the training time interval length.

\begin{figure}[b!]
\begin{center}
    \large
    \textbf{$\boldsymbol{\mu}\ \mathbf{= 2.45}$}\\
    \vspace{1mm}
    \normalsize \textbf{Single species trend}
\end{center}
\vspace{-0.1cm}
\hspace{-0.6cm}
\centerline{
    \includegraphics[width=0.3\textwidth]{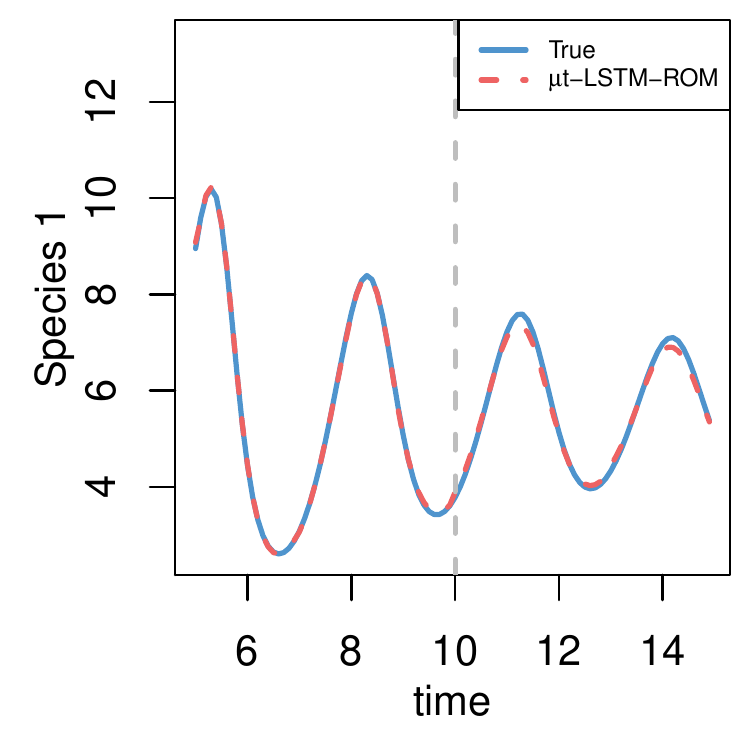}
    \includegraphics[width=0.3\textwidth]{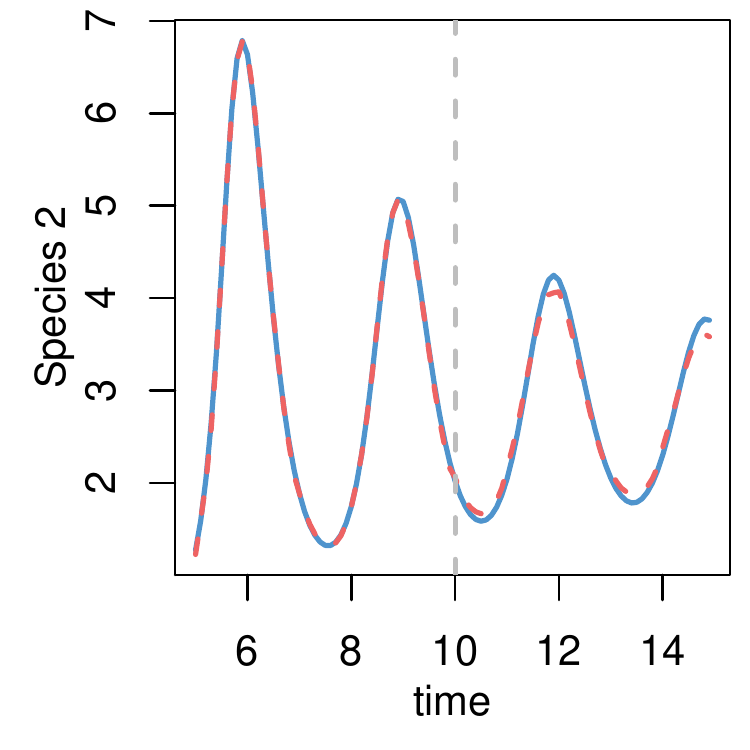}
    \includegraphics[width=0.3\textwidth]{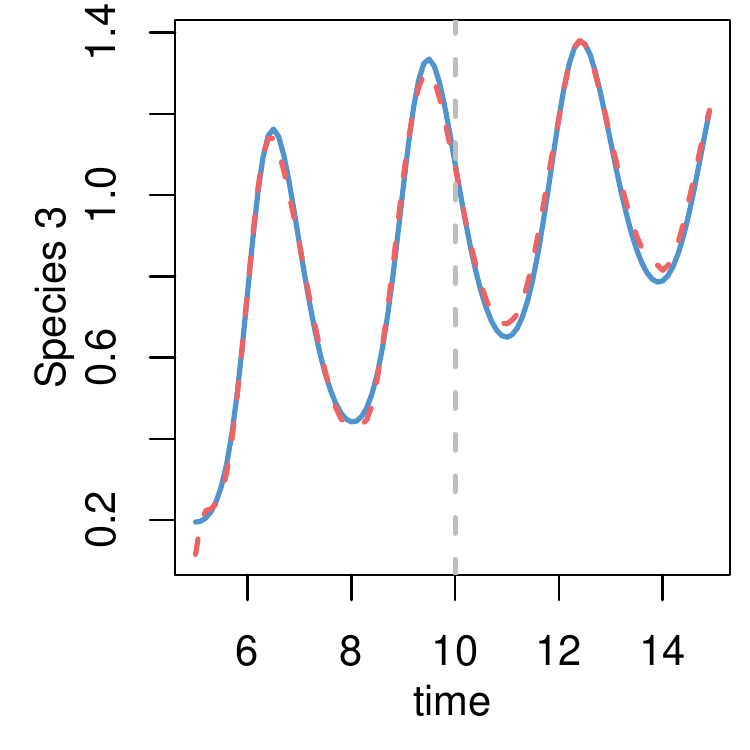}
}
\vspace{-0.25cm}
\begin{center}
    \textbf{Phase space plots}
\end{center}
\vspace{-0.15cm}
\hspace{0.7cm}
\centerline{
\begin{minipage}[]{0.49\textwidth}
    \includegraphics[width=5cm]{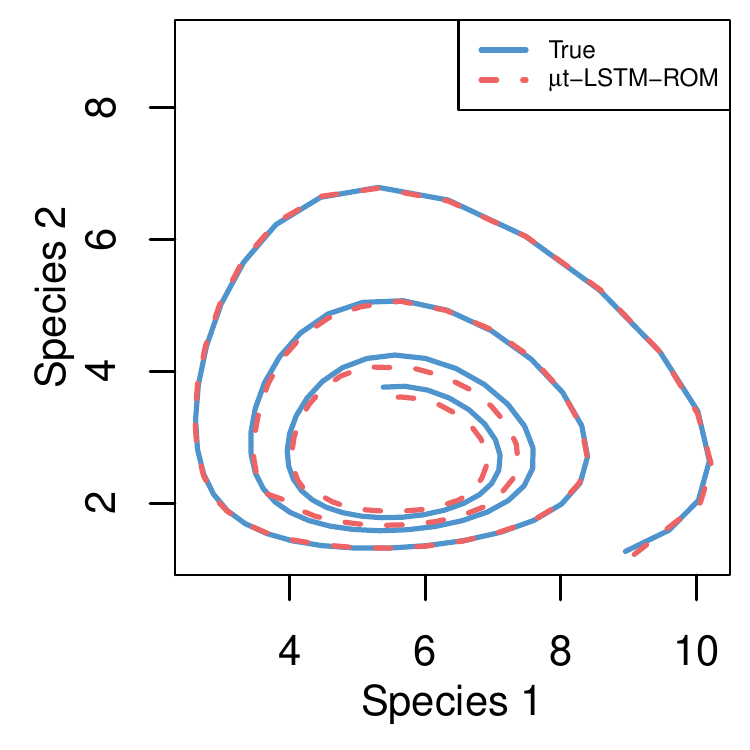}\\
\end{minipage}
\begin{minipage}[]{0.49\textwidth}
    \includegraphics[width=5cm]{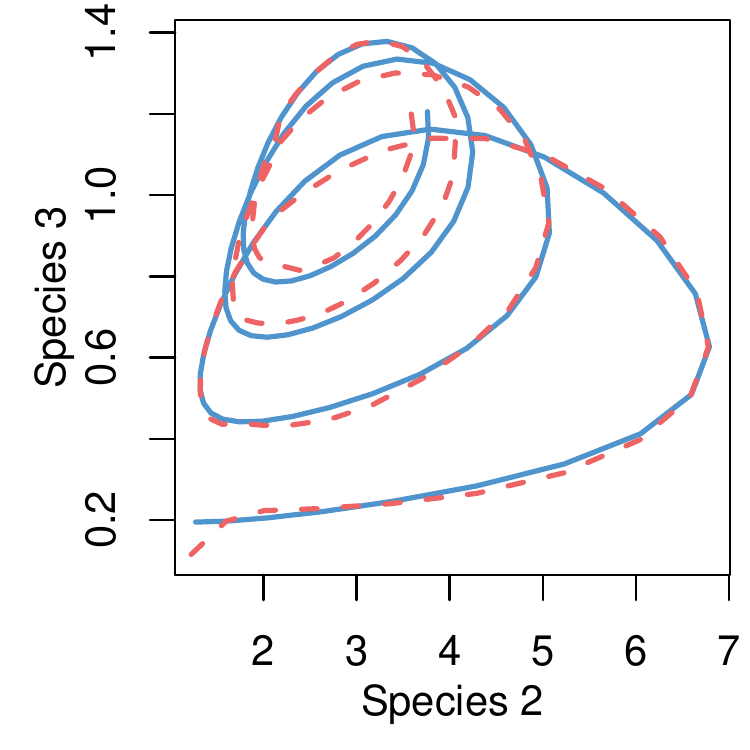}\\
\end{minipage}
}
\vspace{-0.75cm}
\caption{Test case 1 -- Lotka-Volterra system. Simulation results considering the aperiodic Lotka-Volterra model for $\mu\ = 1.25$, $\mu t$-LSTM-ROM framework. Grey dotted line indicates the starting extrapolation time.}
\label{extrapLV}
\end{figure}

The $t$-LSTM-ROM architecture used in addition to the LSTM-ROM one alredy described considers $p=10$ previous time steps in order to make inference on a $k=10$ time steps horizon. Its training, performed without accounting for early stopping, took 1000 epochs (355s). The plots in Figure \ref{extrapLV} show the extrapolation performances of the method when applied to this test case.

\begin{figure}[h!bt]
\centerline{
\includegraphics[width=0.875\textwidth]{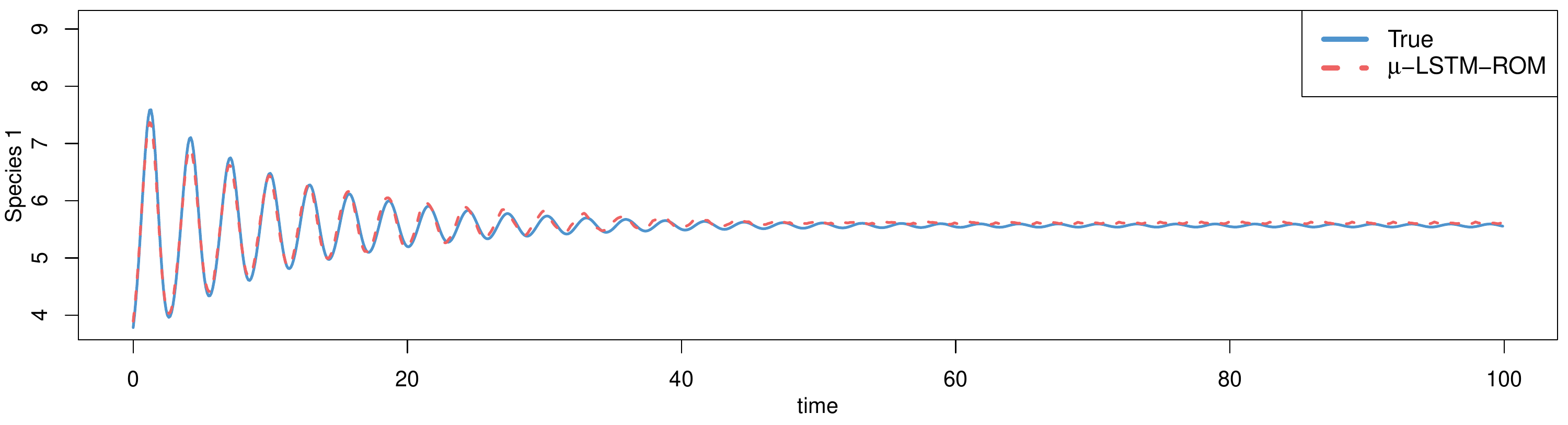}
}
\centerline{
\includegraphics[width=0.875\textwidth]{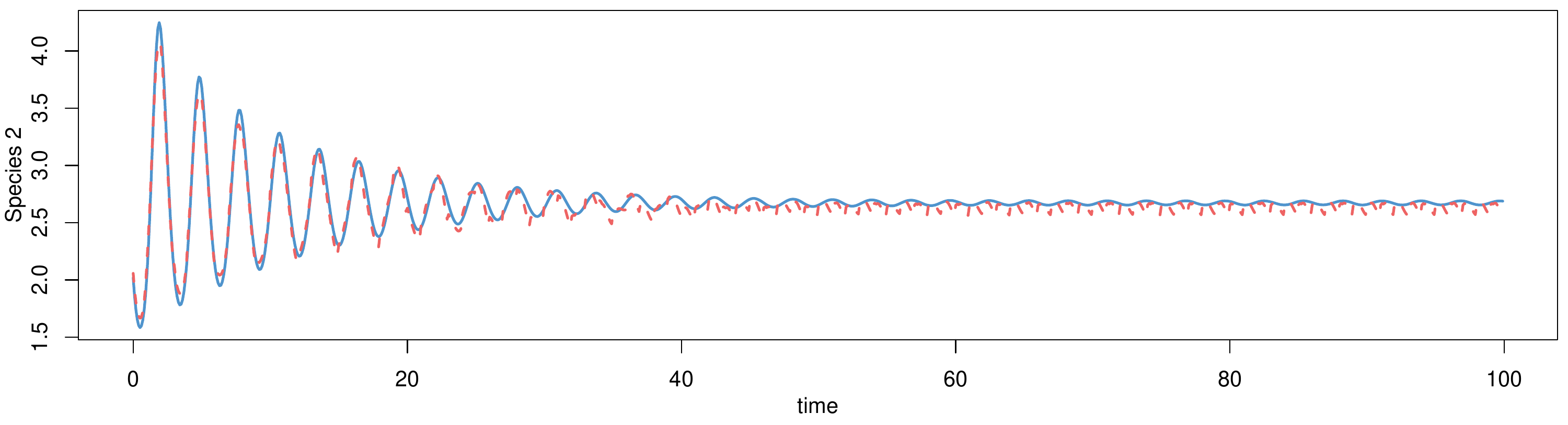}
}
\centerline{
\includegraphics[width=0.875\textwidth]{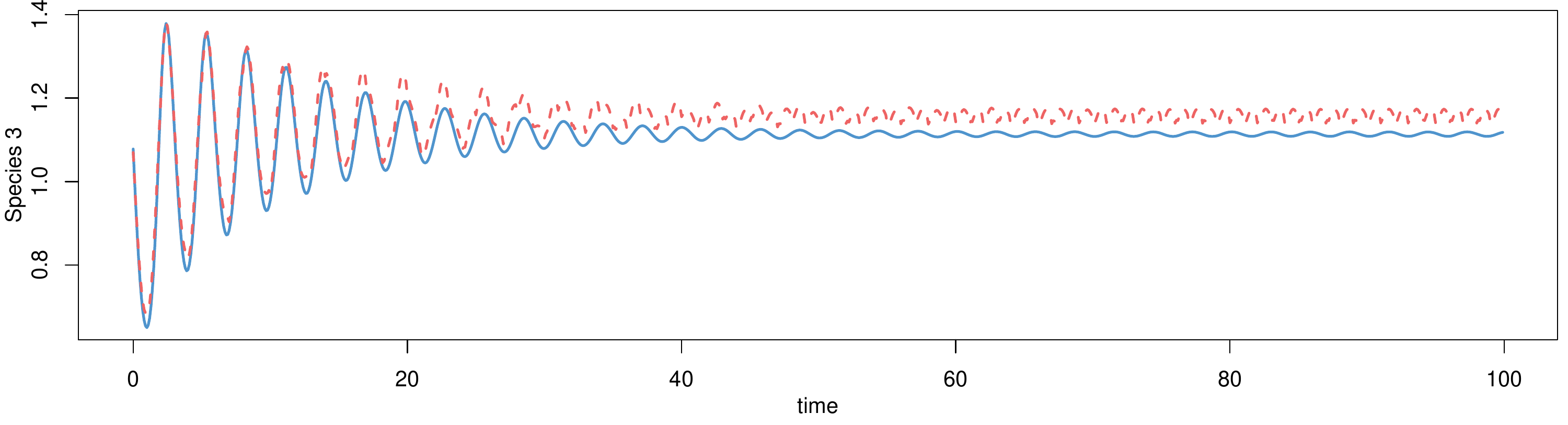}
}
\vspace{-0.3cm}
\caption{Test case 1 -- Lotka-Volterra system. Long term ($T_{ext} = 100$) time extrapolation for the evolution of the three species in the aperiodic Lotka-Volterra case.}
\label{longrunLV}
\end{figure}

For the problem at hands the relative error $\epsilon_k$ obtained in the aforementioned time extrapolation context can be summarized by Table \ref{table3}. Extrapolation performances of the presented framework are very satisfying also considering the long run, as shown for $\mu = 2.45$ in Figure \ref{longrunLV}. In particular, the temporal evolution of the problem on a time window which is 10 times larger than the training domain is predicted with remarkable accuracy as well as the equilibrium asymptotes.

\begin{table}[b!]
\centering
\begin{tabular}{@{}ll@{}}
\toprule
\textbf{$\mathbf{\boldsymbol{\epsilon}_{k}^{mean}}$} & \textbf{$\mathbf{\boldsymbol{\epsilon}_{k}^{max}}$} \\ \midrule
$1.341 \cdot 10^{-3}$                & $3.241 \cdot 10^{-2}$              \\ \bottomrule
\end{tabular}
\caption{Test case 1 -- Lotka-Volterra system. Relative error indicators for the $\mu t$-POD-LSTM-ROM framework applied to the aperiodic Lotka-Volterra system.}
\label{table3}
\end{table}

\subsection{Unsteady advection-diffusion-reaction equation}
\label{ADRsection}

We now consider the case of a parameterized unsteady advection-diffusion-reaction problem. In particular, our goal in this case is to approximate the solution $u=u({\bf x}, t; \boldsymbol{\mu})$ of a linear parabolic PDE initial value problem of the following form:

\begin{equation}
\left\{
\begin{aligned}
& \frac{\partial u}{\partial t}  - \textnormal{div}( \mu_1  \nabla u) + \mathbf{b}(t) \cdot \nabla u + c u   =  f(\mu_2, \mu_3), & \  & (\mathbf{x}, t) \in \Omega \times (0,T),\\
& \mu_1 \nabla u  \cdot \mathbf{n} = 0, &\ & (\mathbf{x}, t) \in \partial \Omega \times (0,T), \\
& u(0) = 0 & \  &  \mathbf{x} \in \Omega.
\end{aligned}
\right.
\label{eq:ADR}
\end{equation}
in the two-dimensional domain $\Omega = (0, 1)^2$, where 
\begin{equation*}
f(\mathbf{x}; \mu_2, \mu_3) = 10 \exp(-((x-\mu_3)^2 + (y-\mu_4)^2)/0.07^2) \qquad \mbox{and} \qquad 
\mathbf{b}(t) = [\cos(t), \sin(t)]^T.
\end{equation*}

Regarding the parameters, we took $n_{\mu} = 3$ parameters belonging to $\mathcal{P} = [0.002, 0.005] \times [0.4, 0.6]^2$. A similar framework has been considered in \cite{fresca2021poddlrom}. In particular we discretized the training parameters space $\mathcal{P}_{train}$ considering $\mu_1 \in \{0.002, 0.003, 0.004, 0.005\}$ and $(\mu_2, \mu_3) \in \{0.40, 0.45, 0.50, 0.55, 0.60\}^2$, for a total $N_{train} = 100$ different instances in $\mathcal{P}_{train}$. Parameters $\mu_2$ and $\mu_3$ influence the location of the distributed source in the spatial domain, while the dependence on $\mu_1$ impacts on the relative importance between the advection and the diffusion terms. 

Note that the dependence of the solution on $\mu_2$ and $\mu_3$ is nonlinear and therefore the problem is nonaffinely parametrized. This would have required the extensive use of hyper-reduction techniques such as, \textit{e.g.}, discrete empirical interpolation method (DEIM) \cite{DEIM1} in order to properly address the construction of a projection-based ROM exploiting, \textit{e.g.} the reduced basis method, thus reducing the performance of this latter.

The FOM snapshots have been obtained by means of a spatial discretization obtained with linear ($\mathbb{P}_1$) finite elements considering $N_h = 10657$ degrees of freedom (DOFs) and a time discretization relying on a Backward Differentiation Formula (BDF) of order 2. The time step used for the time discretization is $\Delta t = 2\pi/20$ on $(0, T)$ with $T=10\pi$.

The $\mu$-POD-LSTM-ROM network used to tackle this problem was trained on $N_t = 100$ time steps for each instance of the parameters space. The dimension of the hidden representation of the LSTM autoencoder has been fixed to $n = 100$. In particular, the network is composed by an initial dense part aimed at reducing the POD dimensions to properly feed an encoder composed by two stacked LSTM cells and then a further dense part is traversed in order to reach the reduced dimension state. The decoder is then symmetrically composed by a dense part, two stacked LSTM cells and another dense network to expand the output up to the required $N$ dimension, for a total number of trainable parameters equals to $|\boldsymbol{\theta}| = 402154$.

The POD dimension was set to be $N = 64$ and the randomized version of the POD has been performed in order to reduce the computational effort required by this stage. After POD, the reduced order snapshots have been scaled in the $[0, 1]$ range. This choice has been taken because of the low magnitude of the solution, ranging in $[0, 0.1]$. In such cases, a \textit{MinMax normalization}, defined as
\begin{equation*}
    x' = \frac{x - min(x)}{max(x) - min(x)} \ \in [0,1]
\end{equation*}
has proven to be useful when applied to the reduced order vectors for maximizing neural network performances \cite{scale}.
Note that, according to our experiments, the performances of the POD-DL-ROM improve significantly when scaling the input, while $\mu$-POD-LSTM-ROM seems to be more robust to non-scaled data.

\begin{figure}[t]
\begin{center}
    \large \textbf{Simulation results - $\boldsymbol{\mu}$ = (0.0025, 0.4250, 0.4250)}
\end{center}
\vspace{-0.2cm}
\centering
\includegraphics[width=0.775\textwidth]{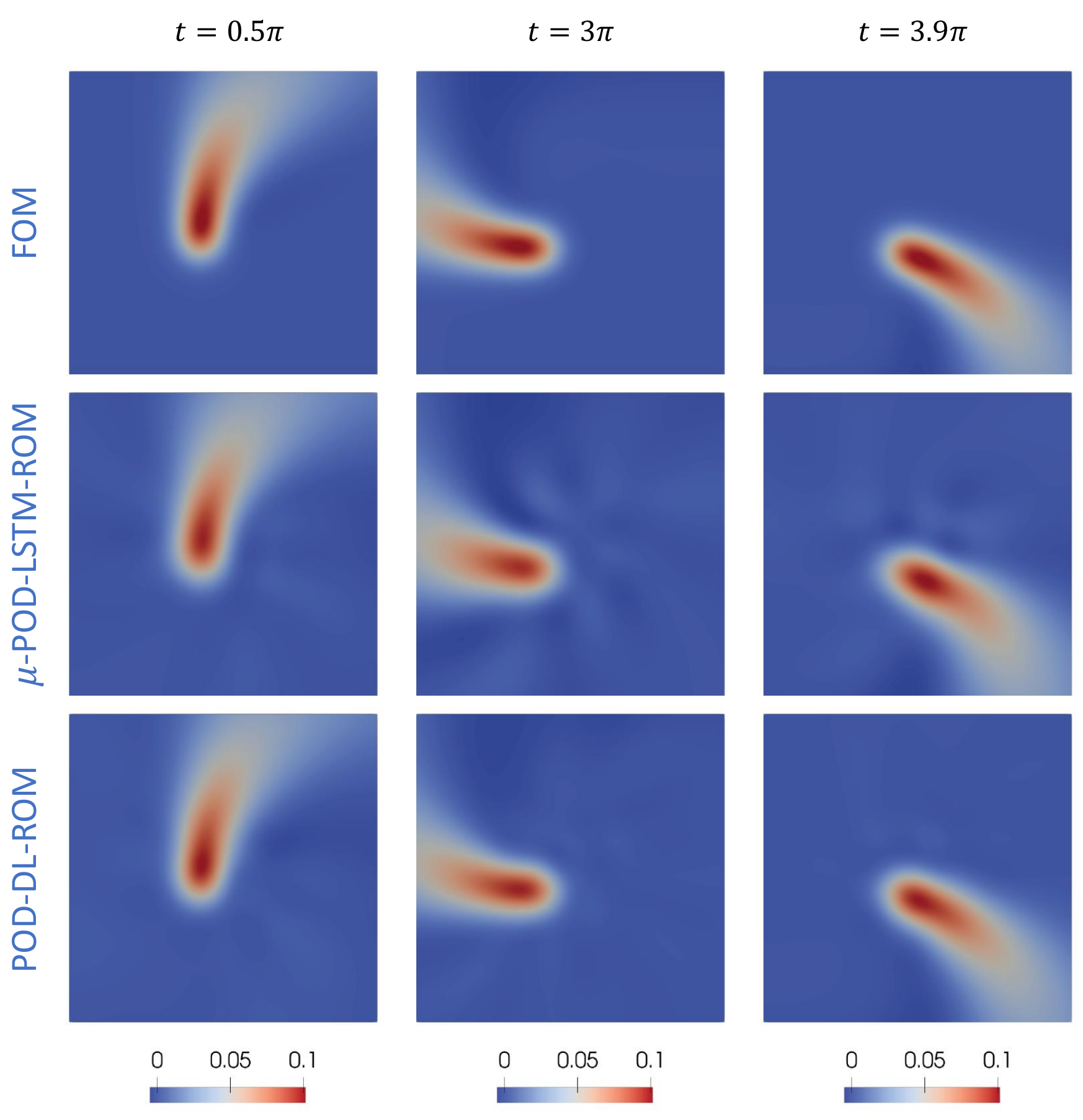}
\caption{Test case 2 -- ADR equation. Simulation results for $\boldsymbol{\mu} = (0.0025, 0.4250, 0.4250)$. Top: FOM results, center: $\mu$-POD-LSTM-ROM results, bottom: POD-DL-ROM results.}
\label{ADRfig1}
\end{figure}

Also in this case, we consider the POD-DL-ROM framework for the comparison with the $\mu$-POD-LSTM-ROM. In particular, the POD-DL-ROM architecture used in this context is the one considered in \cite{fresca2021poddlrom} and it is based on a convolutional architecture aimed at the dimensionality reduction of the input. Therefore, the architecture relies on an initial convolution followed by a feed forward neural network for the encoder, and then symmetrically another feed forward neural network and some deconvolution layers in order to reconstruct the spatial dependence of the solution, for a total number of trainable parameters of $|\boldsymbol{\theta}| = 269207$.
The training of $\mu$-POD-LSTM-ROM took 1243 epochs (2145s), while the one of POD-DL-ROM took 591 epochs (1121s). The obtained POD-DL-ROM and $\mu$-POD-LSTM-ROM results for an instance of the testing parameters space ($\boldsymbol{\mu} = (0.0025, 0.4250, 0.4250)$) are reported in Figure \ref{ADRfig1}.

The error indicator $\epsilon_{rel}^{POD-LSTM-ROM}$ for the LSTM-ROM case is $8.123 \cdot 10^{-2}$, while for the POD-DL-ROM case we find $\epsilon_{rel}^{POD-DL-ROM} = 4.290 \cdot 10^{-2}$, regarding these latter, the obtained results are compatible with those observed in \cite{fresca2021poddlrom}.
Regarding the relative error $\boldsymbol{\epsilon}_k$, we found the results reported in Table \ref{tableerrADR}. The $95\%$ bootstrap confidence intervals for the mean relative error $\epsilon_{k}^{mean}$ are for the $\mu$-POD-LSTM-ROM case $C_{0.95}^{POD-LSTM-ROM} = [2.523 \cdot 10^{-4}, 2.571 \cdot 10^{-4}]$ and $C_{0.95}^{POD-DL-ROM} = [1.303 \cdot 10^{-4}, 1.316 \cdot 10^{-4}]$ in the POD-DL-ROM case, showing a better overall accuracy performance of POD-DL-ROM.

\begin{table}[b!]
\centering
\begin{tabular}{@{}lll@{}}
\cmidrule(l){2-3}
                  & $\mathbf{\boldsymbol{\epsilon}_{k}^{mean}}$ & $\mathbf{\boldsymbol{\epsilon}_{k}^{max}}$ \\ \midrule
\textbf{$\mu$-POD-LSTM-ROM} &  $2.547 \cdot 10^{-4}$                & $3.864 \cdot 10^{-3}$               \\ \midrule
\textbf{POD-DL-ROM}   & $1.309 \cdot 10^{-4}$                 & $2.103 \cdot 10^{-3}$               \\ \bottomrule
\end{tabular}
\vspace{3mm}
\caption{Test case 2 -- ADR equation. Error indicators in comparison between $\mu$-POD-LSTM-ROM and POD-DL-ROM frameworks.}
\label{tableerrADR}
\end{table}

\begin{table}[t!]
\centering
\begin{tabular}{@{}llll@{}}
\cmidrule(l){2-4}
                  & $\mathbf{t_{NN}^{min}}$ & $\mathbf{t_{NN}^{mean}}$ & $\mathbf{t_{NN}^{max}}$ \\ \midrule
\textbf{$\mu$-POD-LSTM-ROM} & 0.0807s               & 0.0892s                & 0.3202s               \\ \midrule
\textbf{POD-DL-ROM}   & 0.1681s                & 0.1857s                 & 0.4973s                \\ \bottomrule
\end{tabular}
\\

\vspace{0.5cm}

\centering
\begin{tabular}{@{}llll@{}}
\cmidrule(l){2-4}
                  & $\mathbf{t_{rec}^{min}}$ & $\mathbf{t_{rec}^{mean}}$ & $\mathbf{t_{rec}^{max}}$ \\ \midrule
\textbf{$\mu$-POD-LSTM-ROM} & 0.7609s               & 0.8437s                & 1.4095s               \\ \midrule
\textbf{POD-DL-ROM}   & 0.8727s                & 0.9864s                 & 3.7111s                \\ \bottomrule
\end{tabular}
\vspace{3mm}
\caption{Test case 2 -- ADR equation. Temporal results for the comparison between $\mu$-POD-LSTM-ROM and POD-DL-ROM.}
\label{table5}
\end{table}

Table \ref{table5} summarizes the testing time performances of the $\mu$-POD-LSTM-ROM in comparison with POD-DL-ROM. Also in this case we report the results in terms both of the neural network involved time and reconstruction time (see the Lotka-Volterra results for details). We obtained a $52.0\%$ decrease in $t_{NN}^{mean}$ and a $14.5\%$ decrease in $t_{rec}^{mean}$ by using $\mu$-POD-LSTM-ROM over POD-DL-ROM. 

\subsubsection*{Advection-Diffusion-Reaction time extrapolation}
Also in this case, we tested the framework to assess its time extrapolation capabilities, both in the short term and in the long term. Such results in this context are remarkable, as traditional ROMs are not capable of time extrapolation for different parameter values than the ones used for training.

Here we consider the same training parameters used before, while reducing the training time domain in order to consider only the first $60$ time steps ($(0, T) = (0, 6\pi)$). A similar procedure is carried out for the test set, that includes snapshots for the same parametric instances $\boldsymbol{\mu} \in \mathcal{P}_{test}$ considered before, but a reduced time domain, chosen in order to contain just the last 60 time steps of the previous one, \textit{i.e.}, $(T_{in}, T_{fin}) = (4\pi, 10\pi)$. In this way, it is possible to test, on unseen parametric instances (in $\mathcal{P}_{test}$), a 40 time steps long time extrapolation.

The $\mu$-POD-LSTM-ROM architecture used to build the $\mu t$-POD-LSTM-ROM framework is the same as the one considered before, while the $t$-POD-LSTM-ROM architecture considers $p=10$ previous time steps in order to predict the following $k=10$ time steps horizon. The number of trainable parameters for this latter network is $|\boldsymbol{\theta}_{ts}| = 151164$; its training took 1000 epochs (156s).

\begin{figure}[t]
\begin{center}
    \large \textbf{Extrapolation results - $\boldsymbol{\mu}$ = (0.0035, 0.4750, 0.4750)}
\end{center}
\vspace{-0.2cm}
\centering
\includegraphics[width=0.775\textwidth]{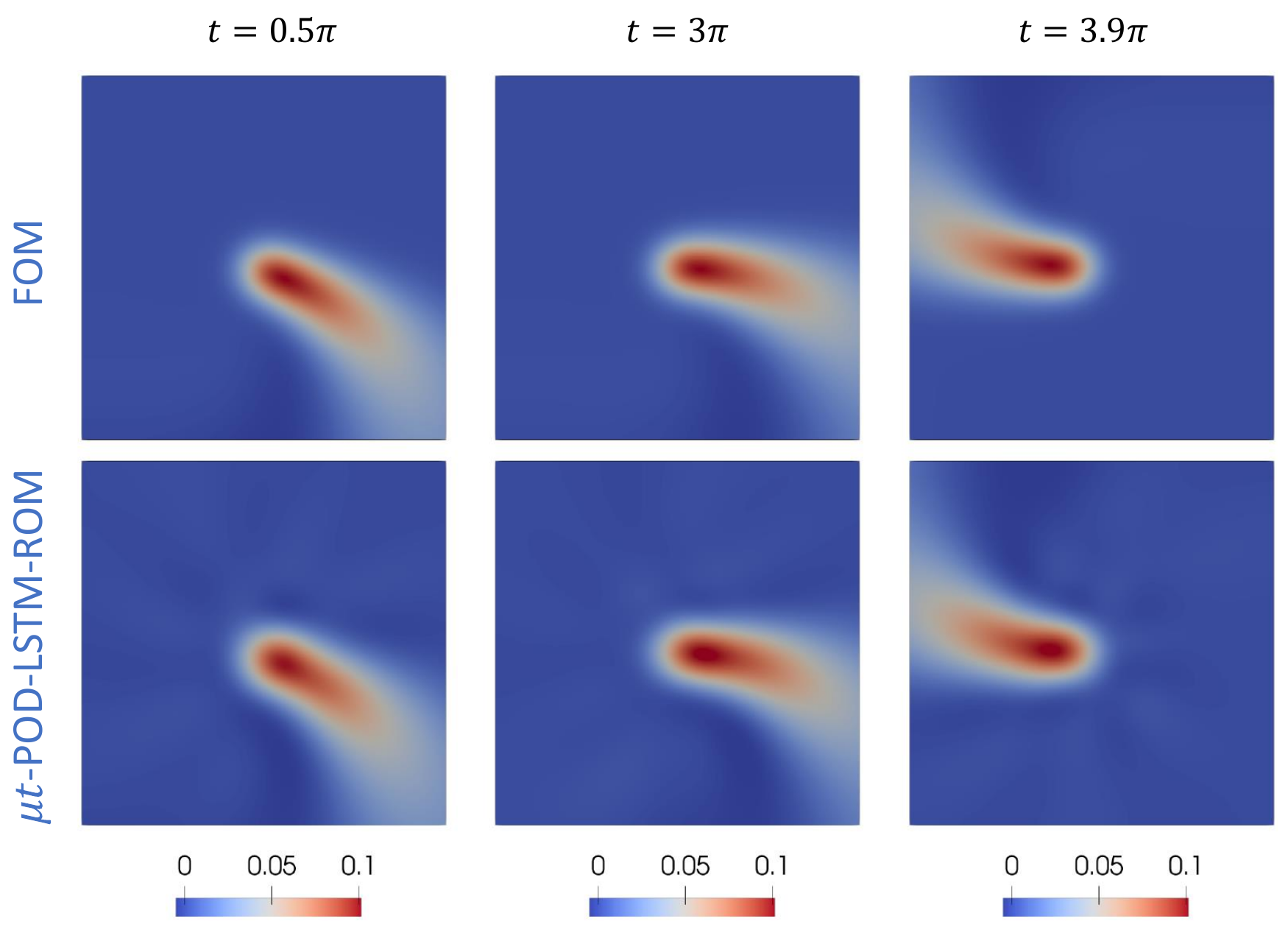}
\caption{Test case 2 -- ADR equation. Simulation results for $\boldsymbol{\mu} = (0.0035, 0.4750, 0.4750)$. Top: FOM results, bottom: $\mu t$-POD-LSTM-ROM results.}
\label{ADRextrap}
\vspace{-0.2cm}
\end{figure}

\begin{figure}[h!bt]
\begin{minipage}[]{0.49\textwidth}
    \begin{center}
    \hspace{1.9cm}
        \footnotesize \textbf{FOM vs $\mu t$-POD-LSTM-ROM}
    \end{center}
        \vspace{0cm}
    \hspace{1.8cm} \includegraphics[width=4.95cm]{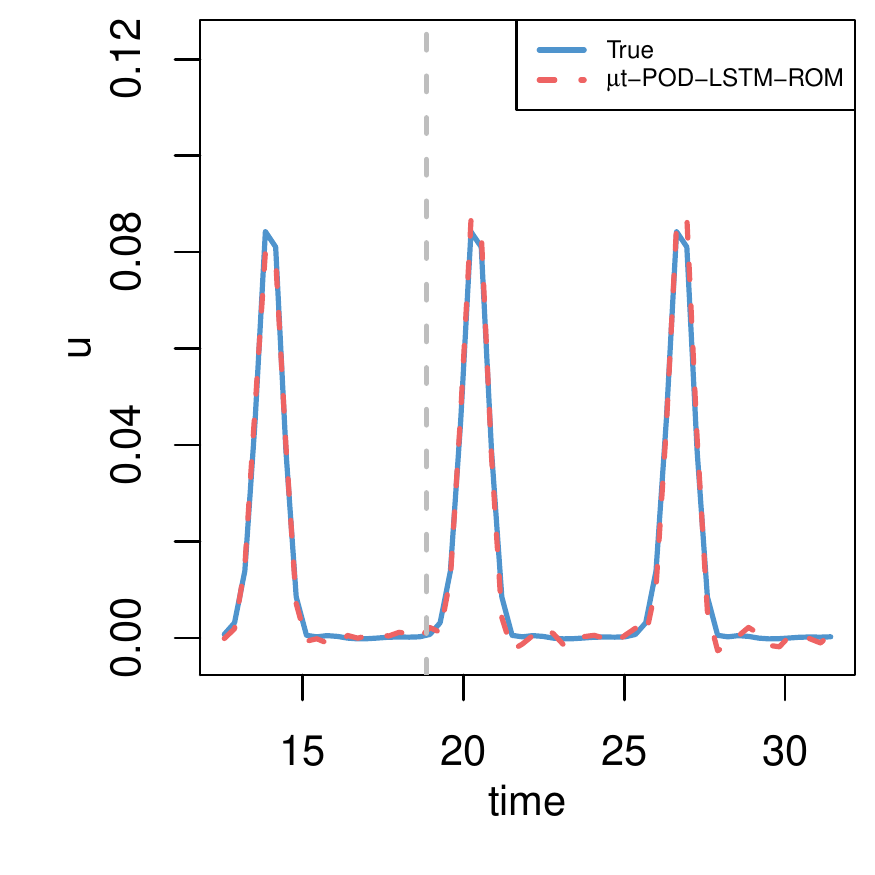}\\
\end{minipage}
\begin{minipage}[]{0.49\textwidth}
    \begin{center}
    \hspace{-1cm}
        \hspace{0.7cm} \footnotesize \textbf{Relative error evolution}
    \end{center}
        \vspace{-0cm}
    \hspace{0.7cm} \includegraphics[width=4.95cm]{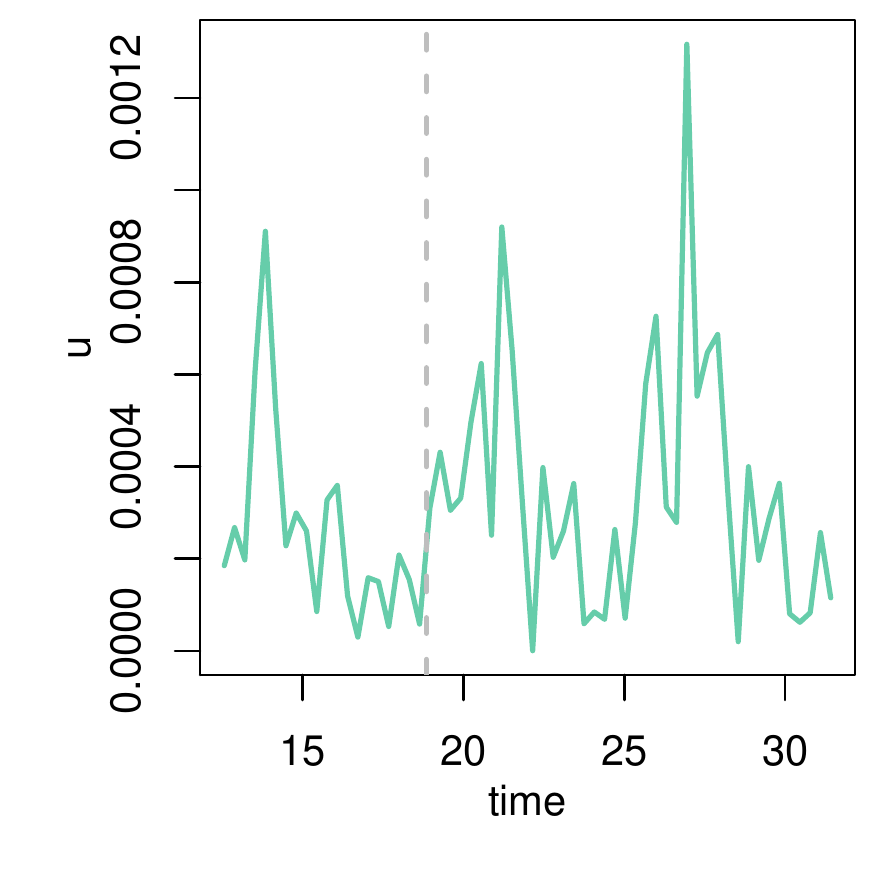}\\
\end{minipage}
\vspace{-0.5cm}
\caption{Test case 2 -- ADR equation. Solution (left) and relative error (right) time evolution when considering advection-diffusion-reaction problem with $\mu t$-POD-LSTM-ROM framework.}
\label{fig:errorADR}
\end{figure}

The obtained results for an instance of the test set ($\boldsymbol{\mu} = (0.0035, 0.4750, 0.4750)$ are reported in Figure \ref{ADRextrap}. Also in this case the extrapolation precision is extremely high, even more considering that traditional POD-based methods do not allow for time extrapolation. Figure \ref{fig:errorADR} reports the time evolution for a single DOF as well as the corresponding relative error evolution, while Table \ref{tab:errorsADR} reports some relevant quantities regarding the relative error for the time extrapolation task.

\begin{table}[b!]
\centering
\begin{tabular}{@{}ll@{}}
\toprule
\textbf{$\mathbf{\boldsymbol{\epsilon}_{k}^{mean}}$} & \textbf{$\mathbf{\boldsymbol{\epsilon}_{k}^{max}}$} \\ \midrule
$2.599 \cdot 10^{-4}$                & $6.261 \cdot 10^{-3}$              \\ \bottomrule
\end{tabular}
\vspace{3mm}
\caption{Test case 2 -- ADR equation. Relative error indicators for the $\mu t$-POD-LSTM-ROM framework applied to the Advection-Diffusion-Reaction problem.}
\label{tab:errorsADR}
\end{table}

\begin{figure}[h!bt]
\centering
\includegraphics[width=0.925\textwidth]{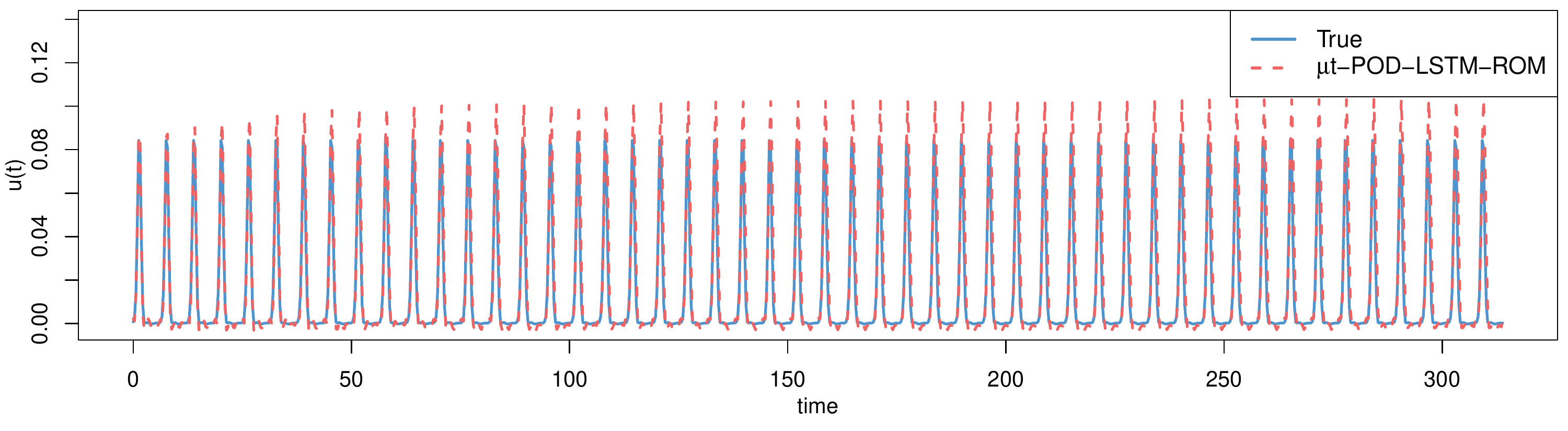}
\vspace{-0.1cm}
\caption{Test case 2 -- ADR equation. Long term (1000 time steps) time extrapolation for the solution at $8430^{th}$ DOF using $\mu t$-POD-LSTM-ROM framework.}
\label{fig:timeextrapADRlong}
\end{figure}

Also in this case, long-term extrapolation capabilities of the framework are extremely satisfying. In Figure \ref{fig:timeextrapADRlong} we report a long-term time extrapolation plot for the temporal evolution of a single DOF that considers a time window 16 times larger than the training domain. In general, over the entire test set performances are satisfying also in the long run, with some issues arising concerning scaling. A maximal systematic error of $\sim 20\%$ arise when such time scales are considered, nonetheless the period is correctly reconstructed and there are no stability issues nor error explosion on the long term. 

\subsection{Unsteady Navier-Stokes equations}
\label{NSsection}

We finally focus on a fluid dynamics example -- the well-known benchmark case of a two-dimensional unsteady flow past a cylinder -- based on incompressible Navier-Stokes equations in the laminar case. Our goal is to approximate the solution $\mathbf{u} = \mathbf{u}(\mathbf{x}, t; \mu)$ of the following problem:
\begin{equation}
\left\{
\begin{aligned}
& \rho \frac{\partial \mathbf{u}}{\partial t} + \rho \mathbf{u} \cdot \nabla \mathbf{u}   - \nabla \cdot \boldsymbol{\sigma}(\mathbf{u}, p) =  \mathbf{0} & \qquad  & (\mathbf{x},t) \in \Omega \times (0,T), \\
&\nabla \cdot \mathbf{u} = 0  & \qquad & (\mathbf{x},t) \in \Omega\times (0,T),\\
& \mathbf{u} =  \mathbf{0} & \qquad & (\mathbf{x},t) \in \Gamma_{D_1} \times (0,T),\\
& \mathbf{u} =  \mathbf{h} & \qquad & (\mathbf{x},t) \in \Gamma_{D_2} \times (0,T),\\
& \boldsymbol{\sigma}(\mathbf{u}, p) \mathbf{n}= \mathbf{0} & \qquad & (\mathbf{x},t) \in \Gamma_{N}\times (0,T), \\
&  \mathbf{u}(0) = \mathbf{0} & \qquad & \mathbf{x} \in \Omega, \ t = 0.
\end{aligned}
\right.
\label{eq:NS}
\end{equation}

\begin{figure}[t]
    \centering
    \includegraphics[width=0.7\textwidth]{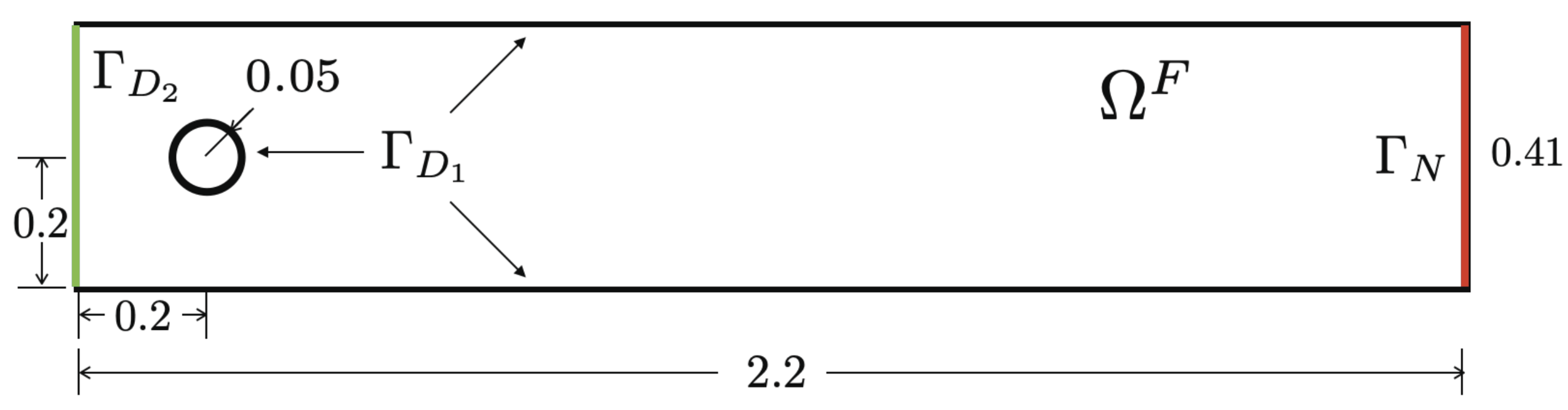}
    \caption{The 2D domain considered for the unsteady Navier-Stokes equations test case.}
    \label{fig:NSdomain}
\end{figure}

The reference domain represents a 2D pipe containing a circular obstacle with radius $r=0.05$ centered in $\mathbf{x_{obs}} = (0.2, 0.2)$, \textit{i.e.}, $\Omega = (0, 2.2) \times (0, 0.41) \char`\\ \bar{B}_r(0.2,0.2)$ (see Figure \ref{fig:NSdomain} for reference); this is a well-known benchmark test case already addressed in \cite{fmfluids1,fresca2021poddlrom}. The domain's boundary is $\partial \Omega = \Gamma_{D_1} \cup \Gamma_{D_2} \cup \Gamma_{N} \cup \partial B_{0.05}(0.2, 0.2)$, where   
$\Gamma_{D_1} = \{ x_1 \in [0, 2.2], x_2 = 0\} \cup \{ x_1 \in [0, 2.2], x_2 = 0.41\}$, $\Gamma_{D_2} = \{ x_1 = 0, x_2 \in [0, 0.41] \}$, and  $\Gamma_N = \{x_1 = 2.2, x_2 \in [0, 0.41] \}$; $\mathbf{n}$ denotes the outward directed versor, normal w.r.t. $\partial \Omega$. We consider $\rho = 1\textnormal{ kg/m}^3$ to be the (constant) fluid density, and denote by
\begin{equation*}
\boldsymbol{\sigma}(\mathbf{u}, p) = -p \mathbf{I} + 2 \nu \boldsymbol{\epsilon}(\mathbf{u})
\end{equation*}
 the stress tensor; here $\nu$ is the fluid's dynamic viscosity, while $\boldsymbol{\epsilon}(\mathbf{u})$ is the strain tensor, \vspace{-0.1cm}
\begin{equation*}
\boldsymbol{\epsilon}(\mathbf{u}) = \frac{1}{2} \big( \nabla \mathbf{u} + \nabla \mathbf{u} ^T \big). \vspace{-0.1cm}
\end{equation*}
We assign no-slip boundary conditions on $\Gamma_1$, while a parabolic inflow profile  \begin{equation}
\label{eq:h}
\mathbf{h}(\mathbf{x},t; \mu) = \left( \frac{4 U(t, \mu) x_2 (0.41-x_2)}{0.41^2} , 0 \right), \qquad \mbox{with } \ \  U(t; \mu) = \mu \sin(\pi t / 8),
\end{equation}
is prescribed at the inlet $\Gamma_{D_2}$; zero-stress Neumann conditions are imposed at the outlet $\Gamma_N$. In this problem, we consider a single parameter ($n_\mu=1$), $\mu \in \mathcal{P}=[1,2]$, which is related with the magnitude of the inflow velocity and directly reflects on the Reynolds number; this latter then varies in the range $\textnormal{Re} \in [66,133]$. Equations (\ref{eq:NS}) have been discretized in space by means of linear-quadratic $(\mathbb{P}_2-\mathbb{P}_1$), inf-sup stable, finite elements, and in time through a BDF of order 2 with semi-implicit treatment of the convective term (see, e.g., \cite{forti2015} for further details) over the time interval $(0,T)$ with $T=6$, considering a time-step $\Delta t = 2 \times 10^{-3}$.

For the sake of training speed we consider $N_t =300$ uniformly distributed time instances and take $N_{train} = 21$ different parameter instances uniformly distributed over $\mathcal{P}$ and $N_{test} = 3$ parameters instances for testing, with  $\mathcal{P}_{test} = \{1.025, 1.725, 1.975\}$, in order to perform testing considering both the center and the boundary of the parameters domain $\mathcal{P}$. In order to assess time extrapolation capabilities of $\mu$-POD-LSTM-ROM, we consider a testing time domain consisting in $N_t = 340$ time steps over the time interval $t \in (0, 6.8]$, resulting in an extrapolation time window of $13.3\%$ w.r.t. the training interval. In this test case, scaling proved to deliver worse results in terms of accuracy and therefore no scaling is performed on POD-reduced data.
We are interested in reconstructing the velocity field, for which  the FOM dimension is equal to $N_h = 32446 \times 2 = 64892$, selecting $N = 256$ as  dimension of the rPOD basis  for each of the two velocity components.

We highlight the possibility, by using $\mu$-POD-LSTM-ROM, to reconstruct only the field of interest, \textit{i.e.}, the velocity $\mathbf{u}$, without the need of taking into account the approximation of the pressure $p$. 

The $\mu$-POD-LSTM-ROM architecture used for this benchmark case consists of a LSTM autoencoder built considering a single LSTM cell both for the encoder and the decoder, without further reducing the dimensionality of the POD-reduced vectors fed to the network. The hidden representation dimension was set to be $n=200$, while the sequence length used for the training was fixed to $K=20$. The number of trainable parameters for the network is $|\boldsymbol{\theta}| = 1365947$.
The $t$-POD-LSTM-ROM architecture used for time extrapolation considers $p=10$ time steps in the past in order to build the inference on $k=10$ time steps in the future, exploiting a single LSTM cell autoencoder and a 3-layers feedforward regressor, for a total number of trainable weights $|\boldsymbol{\theta}_t| = 1616692$. The training of the $\mu$-POD-LSTM-ROM network took 2415 epochs (total time: 5692s), while the training of the $t$-POD-LSTM-ROM architecture took 325 epochs (total time: 694s).

\begin{figure}[t]
    \begin{center}
        \vspace{0.2cm}
    \large \textbf{Extrapolation results}\\
    \normalsize $\mu$ = 1.025
\end{center}
    \centering
    \vspace{-0.35cm}
    \hspace{-0.15cm}
    \includegraphics[width=0.995\textwidth]{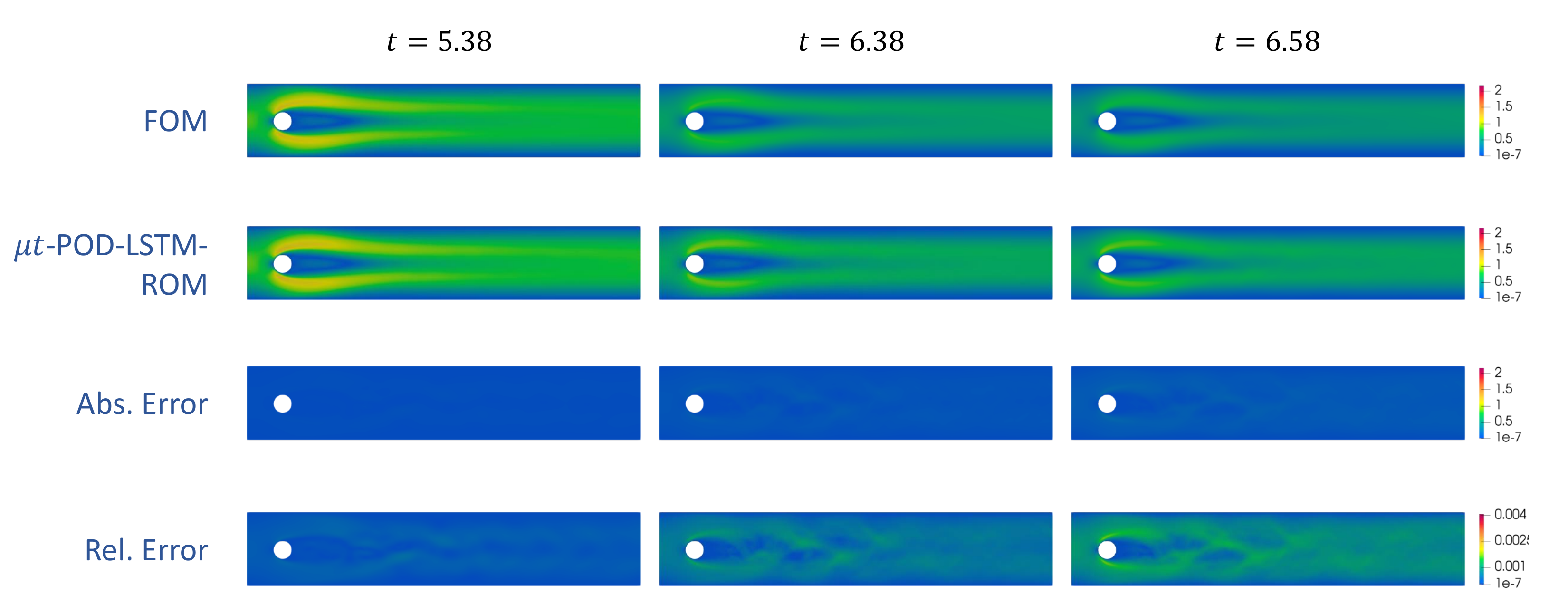}
\begin{center}
    \normalsize $\mu$ = 1.725
\end{center}
    \centering
    \vspace{-0.35cm}
    \hspace{-0.15cm}
    \includegraphics[width=0.995\textwidth]{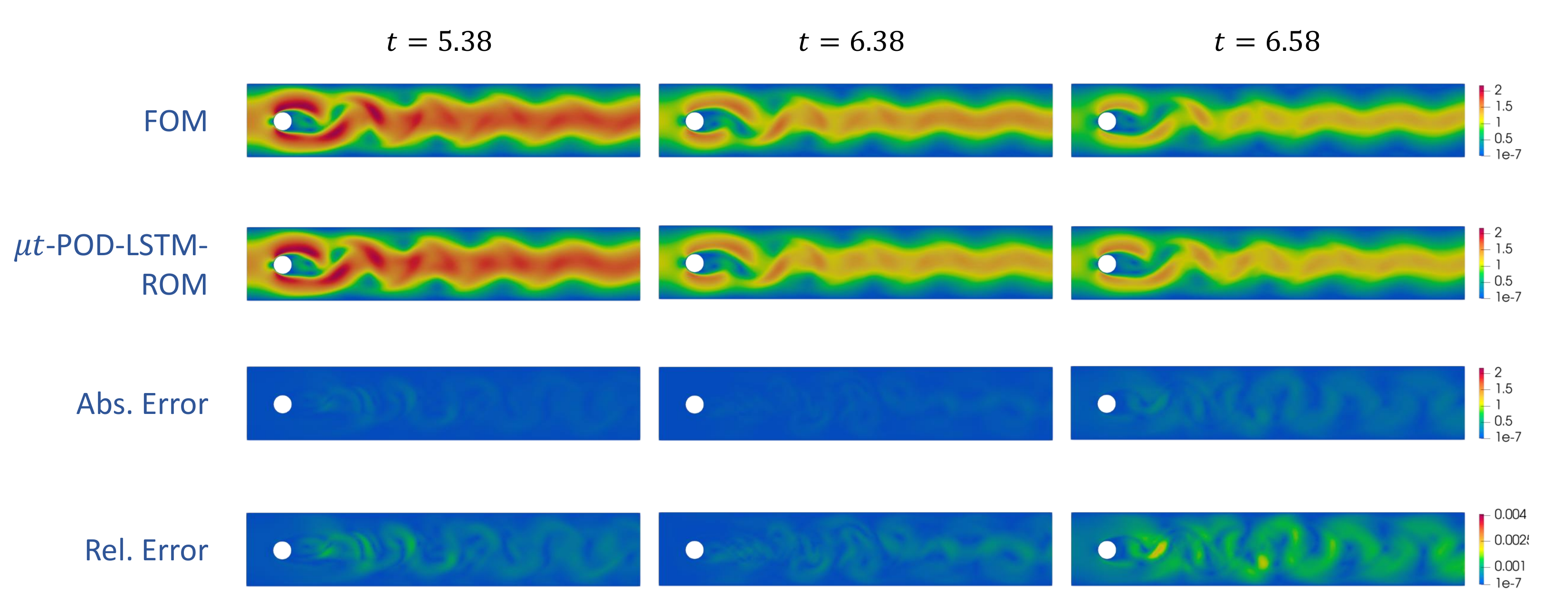}
        \vspace{-0.15cm}
    \caption{Test case 3 -- Navier-Stokes equations. $\mu t$-POD-LSTM-ROM results. Top: $Re = 68$ -- no vortex shedding; Bottom: $Re = 117$ -- vortex shedding. Time extrapolation starts at $t = 6$.}
    \label{fig:NSresults}
\end{figure}

The obtained results for 2 instances of the test set ($\mu = 1.025$, near the border in $\mathcal{P}_{test}$ and no vortex shedding, and $\mu = 1.725$, central in $\mathcal{P}_{test}$ with vortex shedding) are reported in Figure \ref{fig:NSresults} together with absolute and relative errors. A good accuracy is obtained considering time extrapolation, as the relative error reported in Table \ref{tab:NSRelErr} shows. The mean relative error bootstrap 0.95 confidence interval is $C_{0.95}^{POD-LSTM-ROM} = [9.029 \cdot 10^{-5}, 1.053 \cdot 10^{-4}]$. Error indicator for this test case is $\epsilon_{rel} = 5.806 \cdot 10^{-2}$. The little magnitude of the error is especially remarkable considering {\em (a)} the failure of POD-based methods in time extrapolation, {\em (b)} the fact that architecture is just informed on the velocity field and thus cannot take advantage of data on pressure to increase the accuracy of the velocity field prediction and {\em (c)} the complexity of the problem at hand.

\begin{table}[t!]
\centering
\begin{tabular}{@{}ll@{}}
\toprule
\textbf{$\mathbf{\boldsymbol{\epsilon}_{k}^{mean}}$} & \textbf{$\mathbf{\boldsymbol{\epsilon}_{k}^{max}}$} \\ \midrule
$9.808 \cdot 10^{-5}$                & $5.460 \cdot 10^{-3}$              \\ \bottomrule
\end{tabular}
    \vspace{-0.1cm}
\caption{Test case 3 -- Navier-Stokes equations. Relative error indicators for the $\mu t$-POD-LSTM-ROM framework applied to the unsteady Navier-Stokes problem.}
\label{tab:NSRelErr}
\end{table}

\section{Conclusions}
In this work we introduced $\mu t$-POD-LSTM-ROM, a novel non-intrusive LSTM-based ROM framework that extends previous DL-based ROMs with time extrapolation capabilities. In addition, we drastically improved the online performance of the already faster than real-time POD-DL-ROM framework. The strategy followed to pursue our goal splits the solution approximation problem into two parts: {\em (a)} the prediction of the solution for a new parameter instance and {\em (b)} the time extrapolation problem. We therefore introduced two different LSTM-based architectures to address tasks {\em (a)} and {\em (b)} separately. In this way we have been able to replicate the extremely good approximation performances of POD-DL-ROM on unseen parameters instances and -- more importantly -- we enriched the pre-existing DL-based ROMs with time extrapolation capabilities, otherwise hardly obtainable with POD-based ROMs.

We assessed the approximation accuracy, the computational performances and the time extrapolation capabilities on three different test cases: {\em (i)}  a Lotka-Volterra 3 species prey-predator model, {\em (ii)} a linear unsteady advection-diffusion-reaction equation and {\em (iii)} the nonlinear unsteady Navier-Stokes equations with laminar flow. In particular, we observed extremely accurate time extrapolation capabilities, even on a longer term, for both periodic and simple aperiodic cases. Satisfactory time extrapolation capabilities ($\approx 15\%$ of the training time domain) have also been obtained on complex test cases such as the considered benchmark test case in fluid dynamics. An immediate implication of this result is the possibility to produce the FOM snapshots on smaller time domains and therefore to observe a performance improvement in the offline phase. Furthermore, we obtained a $52.0\%$ reduction of prediction times while maintaining the same order of magnitude on the relative error with respect of the (already extremely fast) POD-DL-ROM when applying the framework to advection-diffusion-reaction problems. This allowed us to obtain faster than real-time simulations of physical phenomena occurring in time scales of tenths of a second. Remarkably, the outstanding time performances of our novel framework allow for an "in local" training and testing, thus reducing deployment costs by avoiding the usage of cloud GPU clusters. Ultimately, we provided a novel, fast, accurate and robust framework for the online approximation of parametric time-dependent PDEs, trainable also relying on black-box high-fidelity solvers and applicable to problems of interest in different realms.

\section*{Acknowledgements}
This work has been supported by Fondazione Cariplo, Italy, Grant n. 2019-4608.

\bibliographystyle{elsarticle-num}
\bibliography{biblio}

\vspace{1cm}

\section*{Appendix - Algorithms}
In this appendix, we report in extensive form the algorithms used to train and test the $\mu t$-POD-LSTM-ROM framework. In particular, Algorithms \ref{alg:mu_train} and Algorithm \ref{alg:t_train} describe the procedures used for training the $\mu$-POD-LSTM-ROM and $t$-POD-LSTM-ROM architectures respectively, while Algorithm \ref{alg:mu_test} and Algorithm \ref{alg:t_test} outline the testing stage of the two architectures. Finally, Algorithm \ref{alg:mu_t} specifies the workflow required by $\mu t$-POD-LSTM-ROM framework in order to obtain the fast approximation of a parameterized PDE solution with time extrapolation capabilities.

\begin{algorithm}
\caption{$\mu$-POD-LSTM-ROM training algorithm}
\begin{algorithmic}[1]
\Require Parameter matrix $\mathbf{M} \in \mathbb{R}^{(n_{\boldsymbol{\mu}} + 1) \times N_{train} N_t}$, snapshot matrix $\mathbf{S} \in \mathbb{R}^{N_h \times N_{train} N_t}$, sequence length $K > 1$, validation split $\alpha > 0$, starting learning rate $\eta > 0$, batch size $\textnormal{dim}_{batch}$, maximum number of epochs $N_{ep}$, loss parameter $\omega_h$.
\Ensure  Optimal model parameters $\boldsymbol{\theta}^{*,\mu}= (\boldsymbol{\theta}_{FFNN}^{*,\mu}, \boldsymbol{\theta}_{enc}^{*,\mu}, \boldsymbol{\theta}_{dec}^{*,\mu})$.
\vspace{0.3cm}
\State Compute rPOD basis matrix $\mathbf{V}_N$ \;
\State Compute the POD reduced snapshot matrix $\mathbf{U}_N = [\mathbf{u}_1 | \ldots | \mathbf{u}_{N_{train} N_t}]^T = [\mathbf{V}^T_N \mathbf{u}_{h,1} | \ldots | \mathbf{V}^T_N \mathbf{u}_{h,N_{train} N_t}]^T$
\State Assemble the base tensor $\mathbf{T} \in \mathbb{R}^{N_{train}(N_t-K) \times N \times K}$ with $\mathbf{T}(i,j,k) = (\mathbf{u}_N(t_{\alpha_i} + k\Delta t, \boldsymbol{\mu}_{\beta_i}))_j$ with $\alpha_i = i\textnormal{mod}(N_t - K)$, $\beta_i = \frac{i - \alpha_i}{N_t-K}$ and $(\cdot)_j$ denoting the extraction of the $j^{th}$ component from a vector.
\State Assemble the base parameters tensor $\mathbf{L} \in \mathbb{R}^{N_{train}(N_t-K) \times (n_{\boldsymbol{\mu}}+1) \times K}$ with $\mathbf{L}(i,j,k) = {(\textbf{M}[\alpha_i + k, :])_j}$ with $\alpha_i = i\textnormal{mod}(N_t - K)$, $\beta_i = \frac{i - \alpha_i}{N_t-K}$ and $(\cdot)_j$ denoting the extraction of the $j^{th}$ component from a vector.
\State Randomly shuffle $\mathbf{T}$ and $\mathbf{L}$ \;
\State Randomly sample $\alpha N_{train} (N_t - K)$ indices from $\mathbf{I} = \{0, \dots, N_{train} (N_t - K)-1\}$ and collect them in the vector $val\_idxs$. Build $train\_idxs = I \setminus val\_idxs$
\State Split data in $\mathbf{T} = [\mathbf{T}^{train}, \mathbf{T}^{val}]$ and $\mathbf{L} = [\mathbf{L}^{train}, \mathbf{L}^{val}]$ (with $\mathbf{T}^{train/val} = {\mathbf{T}[train/val\_idxs,:,:]}, \mathbf{L}^{train/val} = {\mathbf{L}[train/val\_idxs,:,:]}$)\;
\State Optionally normalize data in $\mathbf{T}$\;
\State Randomly initialize $\boldsymbol{\theta}^0=(\boldsymbol{\theta}_{FFNN}^0, \boldsymbol{\theta}_{enc}^0, \boldsymbol{\theta}_{dec}^0)$\;
\State $n_{e} = 0$
\While{($\neg$early-stopping \textbf{and} $n_{e} \le N_{ep}$)}
\For{$b = 1 : N_{mb}$}
    \State Sample a minibatch $(\mathbf{T}_{N,K}^{batch}, \mathbf{L}^{batch}) \subseteq (\mathbf{T}_{N,K}^{train}, \mathbf{L}^{train})$\;
    \State $\mathbf{\widetilde{T}}^{batch}_{n,K}(\boldsymbol{\theta}_{enc}^{N_{mb} n_{e} + b}) = \boldsymbol{\lambda}_{n}^{enc}(\mathbf{T}^{batch}_{N,K}; \boldsymbol{\theta}_{enc}^{N_{mb} n_{e} + b})$\;
    \State $\mathbf{T}^{batch}_{n,K}(\boldsymbol{\theta}_{FFNN}^{N_{mb} n_{e} + b}) = \boldsymbol{\phi}_n^{FFNN}(\mathbf{L}^{batch}; \boldsymbol{\theta}_{FFNN}^{N_{mb} n_{e} + b})$\;
    \State $\mathbf{\widetilde{T}}^{batch}_{N,K}(\boldsymbol{\theta}_{FFNN}^{N_{mb} n_{e} + b}, \boldsymbol{\theta}_{dec}^{N_{mb} n_{e} + b}) = \boldsymbol{\lambda}_{N}^{dec}(\mathbf{T}^{batch}_{n,K}(\boldsymbol{\theta}_{FFNN}^{N_{mb} n_{e} + b}); \boldsymbol{\theta}_{dec}^{N_{mb} n_{e} + b})$
    \State Accumulate loss (\ref{eq:loss_encoder})  {on $(\mathbf{T}_{N,K}^{batch}, \mathbf{L}^{batch})$} and compute $\widehat{\nabla}_{\theta} \mathcal{J}$\;
  	\State $\boldsymbol{\theta}^{N_{mb} n_{e} + b + 1} = \textnormal{ADAM}(\eta, \widehat{\nabla}_{\theta} \mathcal{J}, \boldsymbol{\theta}^{N_{mb} n_{e} + b})$\;
  \EndFor
  \State Repeat instructions 13-18 on $(\mathbf{T}^{val}_{N,K}, \mathbf{L}^{val})$ with the updated weights $\boldsymbol{\theta}^{N_{mb} n_{e} + b + 1}$
  \State Accumulate loss (\ref{eq:loss_encoder}) on $(\mathbf{T}^{val}_{N,K}, \mathbf{L}^{val})$ to evaluate early-stopping criterion
  \State $n_{e} = n_{e} + 1$
\EndWhile
\end{algorithmic}
\label{alg:mu_train}
\end{algorithm}


\begin{algorithm}
\caption{$\mu$-POD-LSTM-ROM testing algorithm}
\begin{algorithmic}[1]
\Require Testing parameter matrix $\mathbf{M}^{test} \in \mathbb{R}^{(n_{\boldsymbol{\mu}} + 1) \times N_{test} N_t}$, rPOD basis matrix $\mathbf{V}_N$, optimal model parameters 
$(\boldsymbol{\theta}_{FFNN}^{*,\mu}, \boldsymbol{\theta}_{dec}^{*,\mu})$.
\Ensure ROM approximation matrix $\mathbf{\widetilde{S}}_h \in \mathbb{R}^{N_h \times (N_{test} N_t)}$.
\vspace{0.3cm}
\State Build the reduced testing parameter matrix $\mathbf{M}_{red}^{test} \in \mathbb{R}^{(n_{\boldsymbol{\mu}}+1) \times N_{test}\cdot N_t/K}$ s.t. $\mathbf{M}_{red}^{test}[:,i] = M^{test}[:,K\cdot i]$ $\forall i \in \{0, \dots, N_{test}\cdot N_t/K-1\}$
\State Load $\boldsymbol{\theta}_{FFNN}^{*,\mu}$ and $\boldsymbol{\theta}_{dec}^{*,\mu}$\;
\State $\mathbf{T}_{n,K}(\boldsymbol{\theta}_{FFNN}^{*,\mu}) = \boldsymbol{\phi}_n^{FFNN}(\mathbf{M}^{test}_{red}; \boldsymbol{\theta}_{FFNN}^{*,\mu})$\;
\State $\mathbf{\widetilde{T}}_{N,K}(\boldsymbol{\theta}_{FFNN}^{*,\mu}, \boldsymbol{\theta}_{dec}^{*,\mu}) = \boldsymbol{\lambda}_{N}^{dec}(\mathbf{T}_n(\boldsymbol{\theta}_{FFNN}^{*,\mu}); \boldsymbol{\theta}_{dec}^{*,\mu})$
\State Reshape $\mathbf{\widetilde{T}}_{N,K}(\boldsymbol{\theta}_{FFNN}^{*,\mu}, \boldsymbol{\theta}_{dec}^{*,\mu}) \in \mathbb{R}^{N_h \times N_{test}\cdot N_t/K \times K}$ in $\mathbf{\widetilde{S}}_{N}(\boldsymbol{\theta}_{FFNN}^{*,\mu}, \boldsymbol{\theta}_{dec}^{*,\mu}) \in \mathbb{R}^{N_h \times N_{test} N_t}$
\State $\mathbf{\widetilde{S}}_{h}= \mathbf{V}_N\mathbf{\widetilde{S}}_{N}$
\end{algorithmic}
\label{alg:mu_test}
\end{algorithm}

\begin{algorithm}[h!]
\caption{$t$-POD-DL-ROM training algorithm}
\begin{algorithmic}[1]
\Require Parameter matrix $\mathbf{M} \in \mathbb{R}^{(n_{\boldsymbol{\mu}}) \times N_{train} N_t}$, POD reduced snapshot matrix $\mathbf{U}_{N} \in \mathbb{R}^{N \times N_{train} N_t}$, sequence length $K > 2$, prediction horizon $1 \le k < K$, validation split $\alpha > 0$, starting learning rate $\eta > 0$, batch size $\textnormal{dim}_{batch}$, maximum number of epochs $N_{ep}$, loss parameter $\omega_h$.
\Ensure  Optimal model parameters $\boldsymbol{\theta}^{*,t}= (\boldsymbol{\theta}_{FFNN}^{*,t}, \boldsymbol{\theta}_{enc}^{*,t}, \boldsymbol{\theta}_{dec}^{*,t})$.
\vspace{0.3cm}
\State Assemble the base tensor $\mathbf{T} \in \mathbb{R}^{N_{train}(N_t-K) \times N \times K}$ with $\mathbf{T}(i,j,k) = (\mathbf{u}_N(t_{\alpha_i} + k\Delta t, \boldsymbol{\mu}_{\beta_i}))_j$ with $\alpha_i = i\textnormal{mod}(N_t - K)$, $\beta_i = \frac{i - \alpha_i}{N_t-K}$ and $(\cdot)_j$ denoting the extraction of the $j^{th}$ component from a vector.
\State Assemble the base parameters tensor $\mathbf{L} \in \mathbb{R}^{N_{train}(N_t-K) \times n_{\boldsymbol{\mu}} \times K}$ with $\mathbf{L}(i,j,k) = {(\textbf{M}[\alpha_i + k, :(n_{\boldsymbol{\mu}})])_j}$ with $\alpha_i = i\textnormal{mod}(N_t - K)$, $\beta_i = \frac{i - \alpha_i}{N_t-K}$ and $(\cdot)_j$ denoting the extraction of the $j^{th}$ component from a vector.
\State Randomly shuffle by the first dimension $\mathbf{T}$ and $\mathbf{L}$ \;
\State Randomly sample $\alpha N_{train}(N_t - K)$ indices from $\mathbf{I} = \{0, \dots, N_{train}(N_t - K)-1\}$ and collect them in the vector $val\_idxs$. Build $train\_idxs = I \setminus val\_idxs$
\State Split data in $\mathbf{T} = [\mathbf{T}^{train}, \mathbf{T}^{val}]$ and $\mathbf{L} = [\mathbf{L}^{train}, \mathbf{L}^{val}]$ (with $\mathbf{T}^{train/val} = {\mathbf{T}[train/val\_idxs,:,:]}, \mathbf{L}^{train/val} = {\mathbf{L}[train/val\_idxs,:,:]}$)\;
\State Optionally normalize data in $\mathbf{T}$\;
\State Randomly initialize $\boldsymbol{\theta}^0=(\boldsymbol{\theta}_{FFNN}^0, \boldsymbol{\theta}_{enc}^0, \boldsymbol{\theta}_{dec}^0)$\;
\State $n_{e} = 0$
\While{($\neg$early-stopping \textbf{and} $n_{e} \le N_{ep}$)}
\For{$b = 1 : N_{mb}$}
    \State Sample a minibatch $(\mathbf{T}_{N,K}^{batch}, \mathbf{L}^{batch}) \subseteq (\mathbf{T}_{N,K}^{train}, \mathbf{L}^{train})$\;
    \State Consider $\mathbf{T}_0^{batch} = \mathbf{T}_{N,K}^{train}[:,:,:(K-k)]$ and $\mathbf{T}_1^{batch} = \mathbf{T}_{N,K}^{train}[:,:,(K-k):K)]$
    \State $\mathbf{R}^{batch}_{0}(\boldsymbol{\theta}_{enc}^{N_{mb} n_{e} + b}) = \boldsymbol{\lambda}_{n}^{enc}(\mathbf{T}^{batch}_{0}; \boldsymbol{\theta}_{enc}^{N_{mb} n_{e} + b})$\;
    \State $\mathbf{R}^{batch}_{1}(\boldsymbol{\theta}_{FFNN1}^{N_{mb} n_{e} + b}) = \boldsymbol{\phi}(\mathbf{L}^{batch}; \boldsymbol{\theta}_{FFNN1}^{N_{mb} n_{e} + b})$\;
    \State $\mathbf{H}_n^{batch}(\boldsymbol{\theta}_{FFNN}^{N_{mb} n_{e} + b}, \boldsymbol{\theta}_{enc}^{N_{mb} n_{e} + b}) = \boldsymbol{\phi}'([\mathbf{R}^{batch}_{1}, \mathbf{R}^{batch}_{0}]; \boldsymbol{\theta}_{FFNN2}^{N_{mb} n_{e} + b})$
    \State $\mathbf{\widetilde{T}}^{batch}_{N,k}(\boldsymbol{\theta}_{enc}^{N_{mb} n_{e} + b}, \boldsymbol{\theta}_{FFNN}^{N_{mb} n_{e} + b}, \boldsymbol{\theta}_{dec}^{N_{mb} n_{e} + b}) = \boldsymbol{\lambda}_{N}^{dec}(\mathbf{H}_n^{batch}(\boldsymbol{\theta}_{FFNN}^{N_{mb} n_{e} + b}, \boldsymbol{\theta}_{enc}^{N_{mb} n_{e} + b}); \boldsymbol{\theta}_{dec}^{N_{mb} n_{e} + b})$
    \State Accumulate loss (\ref{eq:loss_lstmts})  {on $(\mathbf{T}_{1}^{batch}, \mathbf{L}^{batch})$} and compute $\widehat{\nabla}_{\theta} \mathcal{J}$\;
  	\State $\boldsymbol{\theta}^{N_{mb} n_{e} + b + 1} = \textnormal{ADAM}(\eta, \widehat{\nabla}_{\theta} \mathcal{J}, \boldsymbol{\theta}^{N_{mb} n_{e} + b})$\;
  \EndFor
  \State Repeat instructions 11-18 on $(\mathbf{T}^{val}_{N,K}, \mathbf{L}^{val})$ with the updated weights $\boldsymbol{\theta}^{N_{mb} n_{e} + b + 1}$
  \State Accumulate loss (\ref{eq:loss_lstmts}) on $(\mathbf{T}^{val}_{N,K}, \mathbf{L}^{val})$ to evaluate early-stopping criterion
  \State $n_{e} = n_{e} + 1$
\EndWhile
\end{algorithmic}
\label{alg:t_train}
\end{algorithm}

\begin{algorithm}[h!]
\caption{$t$-POD-LSTM-ROM testing algorithm}
\begin{algorithmic}[1]
\Require Testing parameter matrix (without times) $\mathbf{M}_{-}^{test} \in \mathbb{R}^{n_{\boldsymbol{\mu}} \times N_{test}}$, rPOD basis matrix $\mathbf{V}_N$, $\mu$-POD-LSTM-ROM reduced approximation matrix optimal model parameters $\mathbf{\widetilde{S}}_N = \mathbf{V}_N^T \mathbf{\widetilde{S}}_h \in \mathbb{R}^{N \times (N_{test} N_t)}$, extrapolation starting point $1 \le t_{ext} \le N_t$, extrapolation length $N_{ext} \ge 1$, prediction horizon $1 \le k < K$ (with $K$ being the sequence length used for training), optimal training parameters
$(\boldsymbol{\theta}_{enc}^{*,t}, \boldsymbol{\theta}_{FFNN}^{*,t}, \boldsymbol{\theta}_{dec}^{*,t})$.
\Ensure ROM extrapolation matrix $\mathbf{\widetilde{E}}_h \in \mathbb{R}^{N_h \times (N_{test} \cdot N_{ext} k)}$.
\vspace{0.3cm}
\State Allocate memory for $\mathbf{\widetilde{E}}_{N,K} \in \mathbb{R}^{N_{test} \times N \times ((N_{ext}-1) k + K)}$
\State Initialize $\mathbf{\widetilde{E}}_{N,K}$ by setting $\mathbf{\widetilde{E}}_{N,K}[:,i,:(K-k)] = {\mathbf{\widetilde{S}}_N[:, (i N_t + t_{ext} - K + k):(i N_t + t_{ext})]}$ ${\forall i \in \{0, \dots, N_{test}-1\}}$

\State Load $\boldsymbol{\theta}_{enc}^{*,t}$, $\boldsymbol{\theta}_{FFNN}^{*,t}$ and $\boldsymbol{\theta}_{dec}^{*,t}$\;
    \State $c = 1$
\For{$j = 1:N_{ext} $}
    \State $\mathbf{R}^{test,0}(\boldsymbol{\theta}_{enc}^{*,t}) = \boldsymbol{\lambda}_{n}^{enc}(\mathbf{\widetilde{E}}_{N,K}[:,:,c:(c+K-k)]; \boldsymbol{\theta}_{enc}^{*,t})$\;
    \State $\mathbf{R}^{test,1}(\boldsymbol{\theta}_{FFNN1}^{*,t}) = \boldsymbol{\phi}(\mathbf{M}_{-}^{test}; \boldsymbol{\theta}_{FFNN1}^{*,t})$\;
    \State $\mathbf{H}_n^{test}(\boldsymbol{\theta}_{FFNN}^{*,t}, \boldsymbol{\theta}_{enc}^{*,t}) = \boldsymbol{\phi}'([\mathbf{R}^{test,1}, \mathbf{R}^{test,0}]; \boldsymbol{\theta}_{FFNN2}^{*,t})$
    \State $\mathbf{\widetilde{E}}_{N,K}(\boldsymbol{\theta}_{enc}^{*,t}, \boldsymbol{\theta}_{FFNN}^{*,t}, \boldsymbol{\theta}_{dec}^{*,t})[:,:,(c+K-k):(c+K)] = \boldsymbol{\lambda}_{N}^{dec}(\mathbf{H}_n^{test}(\boldsymbol{\theta}_{FFNN}^{*,t}, \boldsymbol{\theta}_{enc}^{*,t}); \boldsymbol{\theta}_{dec}^{*,t})$
    \State $c = c + k$
\EndFor
\State Consider the matrix containing extrapolation results only, \textit{i.e.}, $\mathbf{\widetilde{E}}_{N,K}^{ext}(\boldsymbol{\theta}_{enc}^{*,t}, \boldsymbol{\theta}_{FFNN}^{*,t}, \boldsymbol{\theta}_{dec}^{*,t}) = \mathbf{\widetilde{E}}_{N,K}(\boldsymbol{\theta}_{enc}^{*,t}, \boldsymbol{\theta}_{FFNN}^{*,t}, \boldsymbol{\theta}_{dec}^{*,t})[:,:,(K-k):]$
\State Reshape $\mathbf{\widetilde{E}}_{N,K}^{ext}(\boldsymbol{\theta}_{enc}^{*,t}, \boldsymbol{\theta}_{FFNN}^{*,t}, \boldsymbol{\theta}_{dec}^{*,t}) \in \mathbb{R}^{N_{test} \times N \times N_{ext}k}$ in $\mathbf{\widetilde{E}}_{N}(\boldsymbol{\theta}_{FFNN}^{*,t}, \boldsymbol{\theta}_{dec}^{*,t}) \in \mathbb{R}^{N \times (N_{test} \cdot N_{ext} k)}$
\State $\mathbf{\widetilde{E}}_{h}= \mathbf{V}_N\mathbf{\widetilde{E}}_{N}$
\end{algorithmic}
\label{alg:t_test}
\end{algorithm}

\begin{algorithm}[h!]
\caption{$\mu t$-POD-LSTM-ROM training-testing algorithm}
\begin{algorithmic}[1]
\Require The same inputs as Algorithms \ref{alg:mu_train}, \ref{alg:mu_test}, \ref{alg:t_train}, \ref{alg:t_test}.
\Ensure Time extended ROM approximation matrix $\mathbf{\widetilde{S}}_h^{ext} \in \mathbb{R}^{N_h \times N_{test} (t_{ext}+N_{ext}k)}$.
\vspace{0.3cm}
\State Train $\mu$-POD-LSTM-ROM architecture according to Algorithm \ref{alg:mu_train}
\State Train $t$-POD-LSTM-ROM architecture according to Algorithm \ref{alg:t_train}
\State Load $(\boldsymbol{\theta}_{enc}^{*,\mu}, \boldsymbol{\theta}_{FFNN}^{*,\mu}, \boldsymbol{\theta}_{dec}^{*,\mu})$ and $(\boldsymbol{\theta}_{enc}^{*,t}, \boldsymbol{\theta}_{FFNN}^{*,t}, \boldsymbol{\theta}_{dec}^{*,t})$ \;
\State Obtain $\mathbf{\widetilde{S}}_h \in \mathbb{R}^{N_h \times N_{test} N_t}$ by the procedure described in Algorithm \ref{alg:mu_test}
\State Obtain $\mathbf{\widetilde{E}}_h \in \mathbb{R}^{N_h \times (N_{test} \cdot N_{ext}k)}$ by the procedure described in Algorithm \ref{alg:t_test}
\State Assemble the time extended ROM approximation matrix $\mathbf{\widetilde{S}}_h^{ext}$ by concatenating by column part of the previous results as in $\mathbf{\widetilde{S}}_h^{ext} = \mathbf{\widetilde{S}}_h[:,:t_{ext}] \oplus {\mathbf{\widetilde{E}}_h}$
\end{algorithmic}
\label{alg:mu_t}
\end{algorithm}

\end{document}